%
%
\magnification=\magstep 1
\baselineskip=15pt     
\parskip=3pt plus1pt minus.5pt
\overfullrule=0pt
\font\hd=cmbx10 scaled\magstep1
\input amssym.def
\input amssym.tex
\input epsf.tex

\font\tenmsam=msam10
\font\sevenmsam=msam7
\font\fivemsam=msam5
\newfam\msamfam         
\textfont\msamfam\tenmsam
\scriptfont\msamfam\sevenmsam
\scriptscriptfont\msamfam\fivemsam

%

\font\tenmsbm=msbm10
\font\sevenmsbm=msbm7
\font\fivemsbm=msbm5
\newfam\msbmfam         
\textfont\msbmfam\tenmsbm
\scriptfont\msbmfam\sevenmsbm
\scriptscriptfont\msbmfam\fivemsbm

%

\font\teneufm=eufm10
\font\seveneufm=eufm7
\font\fiveeufm=eufm5
\newfam\eufmfam        
\textfont\eufmfam\teneufm
\scriptfont\eufmfam\seveneufm
\scriptscriptfont\eufmfam\fiveeufm

%

%
%
%
%
%
%
%
\newcount\amsfamcount 
\newcount\classcount   
\newcount\positioncount
\newcount\codecount
\newcount\n             
\def\newsymbol#1#2#3#4#5{               
\n="#2                                  
\ifnum\n=1 \amsfamcount=\msamfam\else   
\ifnum\n=2 \amsfamcount=\msbmfam\else   
\ifnum\n=3 \amsfamcount=\eufmfam
\fi\fi\fi
\multiply\amsfamcount by "100           
\classcount="#3                 
\multiply\classcount by "1000           
\positioncount="#4#5            
\codecount=\classcount                  
\advance\codecount by \amsfamcount      
\advance\codecount by \positioncount
\mathchardef#1=\codecount}              
\newcount\famcnt 
\newcount\classcnt   
\newcount\positioncnt
\newcount\codecnt
\def\newmathsymbol#1#2#3#4#5{          
\famcnt=#2                      
\multiply\famcnt by "100        
\classcnt="#3                   
\multiply\classcnt by "1000     
\positioncnt="#4#5              
\codecnt=\classcnt              
\advance\codecnt by \famcnt     
\advance\codecnt by \positioncnt
\mathchardef#1=\codecnt}        
\def\span{{\rm span}_{\bf C}}
\def\rto{\raise.5ex\hbox{$\scriptscriptstyle ---\!\!\!>$}}

\def\Psix{{\bf P^6}}

\def\Pfour{{\bf P^4}}
\def\Pthree{{\bf P^3}}

\def\Ptwo{{\bf P^2}}
\def\Pone{{\bf P}^1}
\def\P{{\bf P}}

\def\A{{\cal A}}

\def\SS{{\cal S}}
\def\X{{\cal X}}
\def\O{{\cal O}}
\def\L{{\cal L}}
\def\PP{{\cal P}}

\def\H{{\cal H}}
\def\G{{\cal G}}

\def\F{{\cal F}}
\def\I{{\cal I}}
\def\M{{\cal M}}

\def\E{{\cal E}}
\def\K{{\cal K}}

\def\HHH{{\bf H}}

\def\EY{\E(1)|_Y}

\def\PEY{\P(\EY)}

\def\OP{\O_{\PEY}(1)}
\def\OP1{\O_{\Pone}}

\def\Pic{\mathop{\rm Pic}}
\def\Spec{\mathop{\rm Spec}}

\def\Im{\mathop{\rm Im}}
\def\im{\mathop{\rm im}}

\def\char{\mathop{\rm char}}

\def\Proj{\mathop{\rm Proj}}

\def\mod{\mathop{\rm mod}}

\def\hom{\mathop{\rm Hom}}

\def\boldz{{\bf Z}}

\def\Gm{{\bf G}_m}

\def\dual#1{{#1}^{\scriptscriptstyle \vee}}
\newsymbol\SEMI226F
\def\rtimes{\mathop{\SEMI}}

\def\exact#1#2#3{0\rightarrow#1\rightarrow#2\rightarrow#3\rightarrow0}
\def\exactm#1#2#3{0\rightarrow#1\rightarrow#2\rightarrow#3\rightarrow1}

\def\mapleft#1{\smash{
  \mathop{\longleftarrow}\limits^{#1}}}
\def\mapright#1{\smash{
  \mathop{\longrightarrow}\limits^{#1}}}

\def\mapdown#1{\Big\downarrow
   \rlap{$\vcenter{\hbox{$\scriptstyle#1$}}$}}

%

%

%

%

%

%
\def\ker{\mathop{\rm ker}}
\def\im{\mathop{\rm im}}
\def\span{{\rm span}_{\bf C}}
\def\deg{\mathop{\rm deg}}
\def\dim{\mathop{\rm dim}}
%

%
\newsymbol\QED1003

\newsymbol\SEMI226F
\def\rtimes{\mathop{\SEMI}}
\font\lineten=line10
\def\arrow#1{\lineten\char'#1}

\def\mapright#1{{\smash{\mathop{\longrightarrow}\limits^{#1}}}}
\def\mapleft#1{\smash{\mathop{\longleftarrow}\limits^{#1}}}
\def\mapdown#1{\downarrow\rlap{$\vcenter{\hbox{$\scriptstyle#1$}}$}}



\def\harrowext#1{\hskip-3.5pt\raise2.3pt\hbox
to#1{\hrulefill}\hskip-3.5pt}
%




%
\def\longnwarrow{\hskip-15pt\lower6pt\hbox{\vbox{\hbox{$\nwarrow$}
\vskip-1.5pt\hbox{\hskip7.95pt\arrow100 }}}}
\def\rto{\raise.5ex\hbox{$\scriptscriptstyle ---\!\!\!>$}}

%

%

%

%
\newsymbol\shortmid2370
%


%

%

%

%


\def\CC{{\bf C}}
\def\RR{{\bf R}}
\def\NN{{\bf N}}
\def\vt{\vartheta}
\def\dim{{\rm dim}\,}
\def\deg{{\rm deg}\,}
\def\gcd{{\rm gcd}\,}
\centerline{\hd  Equations of $(1,d)$-polarized Abelian Surfaces}
\medskip
\centerline{\it Mark Gross\footnote{*}{Supported by NSF grant DMS-9400873}
and Sorin Popescu\footnote{**}{Supported by DFG grant Po-514/1-1}}
\medskip
\centerline{16 February, 1997}
\medskip
{\settabs 3 \columns
\+Department of Mathematics&&
Department of Mathematics\cr
\+Cornell University&&
Columbia University\cr
\+Ithaca, NY 14850&&
New York, NY 10027\cr
\+mgross@math.cornell.edu&&
psorin@math.columbia.edu\cr}
\bigskip
\bigskip
{\hd \S 0. Introduction.}

In this paper, we  study the equations of projectively embedded
abelian surfaces with a polarization of type $(1,d)$. Classical results
say that given an ample line bundle $\L$ on an abelian surface $A$,
the line bundle $\L^{\otimes n}$ is very ample for $n\ge 3$, and furthermore,
in case $n$ is even and $n\ge 4$, the generators of the homogeneous ideal
$I_A$ of the embedding of $A$ via $\L^{\otimes n}$ are all quadratic;
a possible choice for a set of generators of $I_A$
are the Riemann theta relations.

On the other hand,  much less is known about embeddings via
line bundles $\L$ of type $(1,d)$, that is line bundles $\L$ which are not
powers of another line bundle on $A$. It is well-known that if $d\ge 5$,
and $A$ is a general abelian surface, then $\L$ is
very ample, while $\L$ can never be very ample for $d<5$.
However, even if $d\ge 5$,
$\L$ may not be very ample for special abelian surfaces.
We will restrict our attention in what
follows only to the general abelian surface and wish
to know what form the equations take for such a
projectively embedded abelian surface.

A few special cases are well-documented in the literature:
$d=4$, in which case the general surface is a singular octic in $\Pthree$,
cf. [BLvS], and $d=5$ in which case the abelian surface
is described as the zero set of a section of the Horrocks-Mumford
bundle [HM], whereas its homogeneous ideal is generated by 3
(Heisenberg invariant) quintics and  15 sextics (cf. [Ma]).
Also, recent work by Manolache and Schreyer [MS] and by Ranestad
[Ra] provides a description of the equations and  syzygies
in the case  $d=7$.

In this paper, we take this question up for larger values of $d$, and
in particular, we prove the following

\proclaim Theorem. The homogeneous ideal of a general $(1,d)$-polarized
abelian surface is generated by quadrics if $d\ge 10$.

Moreover, in case the embedding considered is with level structure of
canonical type, we can give a  precise symmetric form for these quadrics.

Our approach is as follows: given a line bundle $\L$ of type $(1,d)$
on the surface $A$, we consider the product embedding
$A\times A\subseteq \P^{d-1}\times \P^{d-1}$, with $x$'s as coordinates
on the first $\P^{d-1}$ factor and  $y$'s as coordinates on the second
factor, and construct certain families of
matrices $M$ whose entries are bilinear in the
variables $x_0,\ldots, x_{d-1}$, and $y_0,\ldots, y_{d-1}$,
and which will drop rank on $A\times A$.
Thus setting $(y_0:\ldots:y_{d-1})$ to be some point in $A$, suitable
minors of $M$ will vanish on the surface $A$. Furthermore,
by choosing special values for the parameter $(y_0:\ldots:y_{d-1})$,
one can obtain $M$'s which are anti-symmetric, and hence deduce that
suitable pfaffians of $M$ will vanish on $A$.

These matrices prove to be quite ubiquitous: for $d$ even,
we produce a family of ${d\over 2}\times {d\over 2}$-matrices $M$
which provide equations for elliptic normal curves of degree $d$,
their secant varieties, and abelian surfaces of type $(1,d)$.
However,  in the odd case we produce ${d}\times {d}$-matrices $M$ having
similar features. Finally, we remark here that for $d=5$
these matrices have been first introduced by R. Moore [Mo] in
connection with the Horrocks-Mumford bundle on $\P^4$, and then
later used by several other mathematicians in the same context
(e.g. see [Au], [ADHPR1], [ADHPR2]). Their determinants are exactly
the quintic hypersurfaces in $\P^4$ defined as the wedge product
of two linearly independent sections of the Horrocks-Mumford bundle
(cf. [Mo], [Au]).

The matrices are constructed in \S 2 and some of their properties
are also described in the same chapter. In the remaining part of the
paper, we discuss the structure of the ideal of abelian surfaces
by using  degeneration arguments. Thus in \S 3, we review the most basic
facts about degenerations of abelian surfaces and elliptic curves, and in
\S 4, we construct projectively embedded degenerations
using Stanley-Reisner ideals.  The deepest degenerations
of abelian surfaces we make use of are described by Stanley-Reisner
ideals coming from certain triangulations of the
2-torus. Combinatorics then help us to understand the ideals of these
degenerations. In \S 5 we study basic facts and determine
equations and syzygies for secant varieties of elliptic curves, while
in \S 6 we gather all the previous information to obtain results
on the ideals of general abelian surfaces.

This work was originally motivated by attempts to describe the moduli spaces
$\A^{lev}_{(1,d)}$ of abelian surfaces with $(1,d)$-polarization and
level structure of canonical type. In fact, already in this paper
we obtain a great deal of information; in particular we define rational maps
$$\Theta_d:\A^{lev}_{1,d} \rto X_d$$
for some suitable projective varieties $X_d$, and then show that for
$d\ge 10$, these maps are birational onto their image.

For $d<10$, the situation is definitely more complicated;
the ideals of $(1,d)$-polarized abelian surfaces can never
be generated by quadrics,
and a more careful analysis is required. This
will be carried out in a sequel to this paper [GP], and will include
a detailed analysis of the moduli spaces $\A_{(1,d)}^{lev}$
for small values of $d$. In  particular, we will obtain the following:

\proclaim Theorem. $\A^{lev}_{(1,d)}$ is rational for
$6\le d\le 10$ and $d=12$,  unirational and non-rational for $d=11$
(being birational with the Klein cubic
$\{\sum_{i\in\boldz_5}x_i^2x_{i+1}=0\}\subset\P^4$), while
$\A_{(1,d)}$ is unirational for $d=14, 16, 18$ and $20$.

An important point will be that, in most of the cases covered by the
theorem, pencils of abelian surfaces on Calabi-Yau 3-folds will account for
the projective lines contained in these moduli spaces. The sequel
to this paper will contain many more details on the geometry
of these moduli spaces.

{\it Acknowledgments:} We would like to thank Alf Aure, Robert Connelly,
Wolfram Decker, David Eisenbud, Kristian Ranestad, and
Frank-Olaf Schreyer for many useful
discussions. We are also grateful to Mike Stillman and Dave Bayer for
{\sl Macaulay} [BS], which helped us tremendously to understand the
shape of the equations described in this paper.

\par\bigskip
\noindent
{\hd \S 1. Generalities.}

{\bf \S 1.1. Canonical Theta Functions and the Heisenberg Group.}

We review here some basic facts about abelian varieties. We follow [LB].

Let $A$ be an abelian variety of dimension $g$ over the complex numbers,
$A\cong V/\Lambda$, with $V$ a $g$-dimensional complex vector space
and $\Lambda$ a lattice.
Let $\L$ be an ample line bundle on $A$.
We denote by $H$ the first Chern class of $\L$; $H$ is a {\it polarization}
on the abelian variety. As usual, $H$ can be thought of as a positive-definite
Hermitian form, whose imaginary part, $E:=\Im(H)$, takes integer values
on $\Lambda$.
$\L$ induces a natural
map from $A$ to its dual, $\phi_{\L}:A\rightarrow \hat A$,
given by $x\mapsto t_x^*\L\otimes\L^{-1}$,
where $t_x:A\rightarrow A$ is the morphism given by translation by $x\in A$.
The kernel of $\phi_{\L}$, $K(\L)$,
is always of the form
$K(\L)\cong(\boldz/d_1\boldz\oplus\cdots\oplus\boldz/d_g\boldz)^{\oplus 2}$,
where $d_1|d_2|\cdots|d_g$,
and this only depends on $H$. The ordered $g$-uple
$D=(d_1,\ldots,d_g)$ is then called the
{\it type} of the polarization. We write $\boldz^g/D\boldz^g$ for
$\boldz/d_1\boldz\oplus\cdots\oplus\boldz/d_g\boldz$.
On $K(\L)$, the Weil pairing
induced by $E$ is
 given by
$$e^{\L}(x,y)=\exp(2\pi i E(x,y))$$
for $x,y\in K(\L)$. This pairing  depends only on the
polarization $H$ and thus will be denoted in the sequel
mainly as $e^H$.

A decomposition
$$\Lambda=\Lambda_1\oplus\Lambda_2$$
is said to be a {\it decomposition for $\L$} if $\Lambda_1$ and $\Lambda_2$
are isotropic for $E$. This induces a decomposition of real vector
spaces $V=V_1\oplus V_2$. Since $K(\L)=\Lambda(\L)/\Lambda$, where
$$\Lambda(\L)=\{v\in V| E(v,\Lambda)\subseteq\boldz\},$$
a decomposition of $\Lambda$ induces a decomposition
$$K(\L)=K_1(\L)\oplus K_2(\L),$$
with $K_1(\L)\cong K_2(\L) \cong \boldz^g/D\boldz^g$, both subgroups being
isotropic with respect to the Weil pairing.

Given a decomposition, we can define a semicharacter $\chi_0:V\rightarrow
{\bf C}_1$ by
$$\chi_0(v)=\exp({\pi i E(v_1,v_2)})$$
where $v=v_1+v_2$ with $v_i\in V_i$. Via the Appell-Humbert theorem ([LB],
Theorem 2.2.3) this determines a line bundle
$\L_0:=\L(H,\chi_0)$. Then the line bundle $\L$ can be written as $t_c^*\L_0$
for some $c\in V$, unique up to translation by elements of $\Lambda(\L)$.
$c$ is called a {\it characteristic} for the line bundle $\L$,
and $\L_0$ is said to be a line bundle of characteristic zero with respect
to the given decomposition.

Given an ample line bundle $\L$, a decomposition for $\L$, and a choice
of a characteristic $c$ for $\L$, there is a unique basis
$$\{\vt^c_x \mid x\in K_1(\L)\}$$
of {\it canonical theta functions} of the space $\Gamma(X,\L)$.
(See [LB], \S 3.2.) We will often
omit mention of $c$, and write $\vt^{\L}_x$ when we are dealing with several
line bundles at the same time.

Given an element $x\in K(\L)$, we have an isomorphism $t_x^*\L
\cong \L$. In general $x$ then induces a projective automorphism
on $\P(H^0(\L))$ which yields a representation $K(\L)\rightarrow PGL(H^0(\L))$.
This representation does not lift to a linear representation of $K(\L)$,
but it does after taking a central extension of $K(\L)$,
$$1\rightarrow{{\bf C}^*}\rightarrow{\G(\L)}\rightarrow{K(\L)}\rightarrow0,$$
whose Schur commutator map is the previously defined pairing $e^H=e^{\L}$.
$\G(\L)$ is called the {\it theta group} of $\L$. The theta group is isomorphic
to the Heisenberg group $\H(D)$, which can be described as follows:
as a set it is ${\bf C}^*\times K(D)$, where $K(D)\cong
\boldz^g/D\boldz^g\oplus\boldz^g/D\boldz^g$. Let $f_1,\ldots,f_{2g}$
be the standard basis of $K(D)$, and define an alternating multiplicative
form $e^D:K(D)\times K(D)\rightarrow {\bf C}^*$ by
$$e^{D}(f_{\nu},f_{\mu})
:=\cases{\exp(-2\pi i/d_{\nu})&if $\mu=g+\nu$\cr
\exp(2\pi i/d_{\nu})&if $\nu=g+\mu$\cr
1&otherwise.\cr}$$
To define the group structure on $\H(D)$, we take for any
$(\alpha,x_1,x_2),(\beta,y_1,y_2)\in\H(D)$
$$(\alpha,x_1,x_2)(\beta,y_1,y_2):=(\alpha\beta
e^D(x_1,y_2),x_1+y_1,x_2+y_2).$$
An isomorphism between $\G(\L)$ and $\H(D)$ which restricts to the
identity on ${\bf C}^*$ is called a {\it theta structure} for $\L$.
Any such isomorphism induces a symplectic isomorphism between $K(\L)$ and
$K(D)$, that is which preserves the alternating pairings $e^{H}$ and
$e^D$, respectively. Classically, a symplectic
isomorphism $\bar{b}:K(\L)\to K(D)$ is called a
{\it level $D$-structure} on $(A,H)$, or a {\it level structure
of canonical type}. Since we are not considering other kinds
of level structures in this paper, we will refer in the sequel
to these level structures simply as {\it level structures}.

As mentioned above, the theta group has a natural representation
$\G(\L)\rightarrow GL(H^0(\L))$, which lifts uniquely the representation
$K(\L)\rightarrow PGL(H^0(\L))$. Given a choice of theta structure,
this representation is isomorphic to the Schr\"odinger representation of
$\H(D)$, defined as follows. Let $W={\bf C}(\boldz^g/D\boldz^g)$ be
the vector space of complex-valued functions on the set $\boldz^g/D\boldz^g$.
The Schr\"odinger representation $\rho:\H(D)\rightarrow GL(W)$
is given by
$$\rho(\alpha,x_1,x_2)(\gamma)=\alpha e^D(\cdot,x_2)\gamma(\cdot+x_1).$$
This representation is irreducible, and the center ${\bf C}^*$ acts
by scalar multiplication, so it yields a projective representation of
$K(D)$.

More explicitly, for a surface, the Schr\"odinger representation
takes the following form on projective space. Let $D=(d_1,d_2)$.
We can write a basis $\{\delta_{\gamma}\mid\gamma\in\boldz^2/D\boldz^2\}$ of
$W$, where $\delta_{\gamma}$ is the delta function
$$\delta_{\gamma}(\gamma')=\cases{0& if $\gamma\not=\gamma'$;\cr
1&if $\gamma=\gamma'$.\cr}$$
We denote by $\HHH_{d_1,d_2}$
the subgroup of $\H(D)$ generated
by $\sigma_1=(1,1,0,0,0)$, $\sigma_2=(1,0,1,0,0)$, $\tau_1=(1,0,0,1,0)$
and $\tau_2=(1,0,0,0,1)$, and these act on $W$ via
$$\eqalign{
\sigma_1(\delta_{i,j})=\delta_{i-1,j},\quad \sigma_2(\delta_{i,j})=
\delta_{i,j-1},\cr
\tau_1(\delta_{i,j})=\xi_1^{-i}\delta_{i,j},\quad
\sigma_2(\delta_{i,j})=\xi_2^{-j}\delta_{i,j},\cr}
$$
where $\xi_k:=\exp(2\pi i/d_k)$.

In the case that $d_1=1$, both $\sigma_1$ and $\tau_1$ are just the identity,
and we shall denote by $\sigma$ and $\tau$ the generators
$\sigma_2$ and $\tau_2$, and leave off the first index on the variables.

Given a decomposition for $\L$ inducing $K(\L)=K_1(\L)\oplus K_2(\L)$,
a basis of canonical theta functions $\{\vt_{\gamma}\mid \gamma
\in K_1(\L)\}$ for $H^0(\L)$
yields the identification of $H^0(\L)$ and $W$ via $\vt_{\gamma}\mapsto
\delta_{\gamma}$
such that the representations
$\G(\L)\rightarrow GL(H^0(\L))$ and $\H(D)\rightarrow GL(V)$ coincide.
Thus if we map $A$ into $\P(\dual{H^0(\L)})$ using as coordinates
$x_{\gamma}=\vt_{\gamma}$,  $\gamma\in
\boldz^g/D\boldz^g$, the image of $A$ will be invariant under the action of the
Heisenberg group via the Schr\"odinger representation. In particular, if
$A$ is embedded this way in $\P(\dual{H^0(\L)})$, then $H^0(\I_A(n))$ is also a
representation of the Heisenberg group. Moreover, in case $d_1=1$,
this is a representation of weight $n$ (i.e., a central element
$z\in{\bf C}^{\ast}$ acts by multiplication with $z^n$),
and hence all its irreducible components will
have dimension $\ge d_2/\gcd(d_2,n)$ (e.g. see [La], [Mu1]).
These are basic facts which will be used over and over again
in the sequel.

Again, if $d_1=1$, the action of the Heisenberg group $\HHH_d:=\HHH_{1,d}$
on the coordinates of $\P(\dual{H^0(\L)})$ is
$$\eqalign{\sigma(x_i)&=x_{i-1},\cr
\tau(x_i)&=\xi^{-i}x_i.\cr}$$

If one considers $K(\L)$ as a subgroup of the automorphism group of $A$ via
translations, then the order 2 subgroup
$\langle(-1_A)\rangle$ acts on $K(\L)$ by inner automorphisms. We may define
$K^e(\L)$ as $K(\L)\rtimes\langle(-1_A)\rangle$, and then define the extended
theta group $\G(\L)^e$ to be a central extension of $K(\L)^e$ by ${\bf C}^*$.
In fact $\G(\L)^e:=\G(\L)\rtimes\langle(-1_\L)\rangle$. Similarly, one
can introduce an extended Heisenberg group, defined by
$$\H^e(D):= \H(D)\rtimes \langle \iota \rangle,$$
where $\iota$ acts on $\H(D)$ via $\iota(\alpha,x_1,x_2)=
(\alpha,-x_1,-x_2)$. An {\it extended} theta structure is an isomorphism
between $\G^e(\L)$ and $\H^e(D)$ inducing the identity on ${\bf C}^*$.
Each extended theta structure restricts to a theta structure, but
a theta structure does not always come from an extended
theta structure. In fact, a theta structure $b:\G(\L)\rightarrow
\H(D)$ can be extended to an extended theta structure if and only
if it is a {\it symmetric theta structure}, that is
if the diagram
$$\matrix{
\G(\L)&\mapright{(-1)_\L}&\G(\L)\cr
\mapdown{b}&&\mapdown{b}\cr
\H(D)&\mapright{\iota}&\H(D)\cr}$$
commutes. In order for a symmetric theta structure to exist,
$\L$ must be a symmetric line bundle, that is $(-1_A)^*\L\cong \L$.
By Theorem 6.9.5 of [LB], there always exist a finite number of symmetric
line bundles of a given polarization,
each admitting a finite number of symmetric theta structures.

The Schr\"odinger representation $\rho$ of $\H(D)$ extends to a representation
$\rho^e$
of $\H^e(D)$, with $\rho^e(\iota)\in SL^{\pm}(W)=\{M\in GL(W)|\det M=\pm 1\}$.
In the case that $A$ is a surface, $\iota$ acts on $W$ by $\iota(\delta_{i,j})
=\delta_{-i,-j}$. We denote by $\HHH^e_{d_1,d_2}$ the subgroup of
$\H^e(D)$ generated by $\HHH_{d_1,d_2}$ and $\iota$.
In fact, the following holds:
$$\HHH^e_{d_1,d_2}=\HHH_{d_1,d_2}\rtimes \langle \iota \rangle.$$
$\iota$, acting as an involution on $W$, has two eigenspaces,
with eigenvalues $\pm 1$. We will refer to the projectivization
of the positive eigenspace as $\P^+\subseteq \P(W)$, and the negative
eigenspace as $\P^-\subseteq \P(W)$.

In particular, if $D=(1,d)$, then $\P^+$ is given by the equations
$$\{x_i=x_{-i}\mid i\in\boldz/d\boldz\}$$
while $\P^-$ is given by the equations
$$\{x_i=-x_{-i}\mid i\in\boldz/d\boldz\}.$$

{\bf \S 1.2. Linear systems on abelian surfaces.}

We recall in this paragraph general results about linear systems on abelian
surfaces.

Let $A$ be an abelian surface and let $\L$ be an ample line bundle on $A$ of
type $D=(d_1,d_2)$. Then Riemann-Roch gives $h^0(\L)=
{1\over 2}{c_1(\L)}^2=d_1d_2$, and
so $\L$ defines a rational map $\psi_{\L}:A\rto \P^{d_1d_2-1}$.

According to Lefschetz's theorem [Mu2], the line bundle $\L$ is
very ample if $d_1\ge 3$. Moreover, the embedding $\psi_\L$
is projectively normal in this case [Ko].

If $d_1=2$, then there exists an ample line bundle $\M$ on $A$ with $\L=\M^2$,
and moreover $|\L|$ is at least base point free. If $\M$ splits, that is
$(A,\M)\cong (E_1\times E_2, \M_1\boxtimes \M_2)$, where $E_1$ and
$E_2$ are elliptic curves
and $\M_1$ a principal polarization on $E_1$, then $\psi_\L$ is the composition
$$E_1\times E_2\mapright{\psi_{\M_1^2\times \M_2^2}}\Pone\times\P^{d_2-1}
\mapright{\rm Segre}\P^{2d_2-1}.$$
In general $\M$ doesn't split and one distinguishes two cases:
either $d_2=2$ and
$\psi_\L:A\to K\subset\P^3$ is of degree 2 on its image, the
Kummer surface associated to $\M$, or $d_2>2$ and then $\psi_\L$ is an
embedding [Oh], [LN]. Furthermore, in this last case the
embedding is not necessarily projectively normal; see [Oh] for an
explicit necessary and sufficient criterion.
In all the above cases, Mumford's theory of the multiplying-sections-map
shows that in case $\psi_\L$ is an embedding, the homogeneous ideal of
the image is generated by quadrics when $d_1\ge 4$ [Mu1], [Ke], or by quadrics
and cubics for $d_1=3$ [Ke], or by quadrics, cubics and
quartics for $d_1=2$ [Ke]. Moreover, if $2|d_1$ and $d_1\ge 4$ then
Riemann's quadratic theta relations [Mu1], [LB] describe
the image $\psi_\L(A)\subset \P^{d_1d_2-1}$ completely. Similarly,
when $3|d_1$, the cubic theta relations  give a complete set of
equations for the image $\psi_\L(A)\subset\P^{d_1d_2-1}$
(see [LB, Theorem 7.6.6]).

Finally, let $\L$ be of type $(1,d)$. Then the
Decomposition theorem [LB, Theorem 4.3.1]
says that $|\L|$ has base curves if and only if
$(A,\L)\cong (E_1\times E_2, \L_1\boxtimes \L_2)$,
where $\L_1$ is a principal polarization on $E_1$ and
$\L_2$ is of type $(d)$ on $E_2$. Assume from now on that $(A,\L)$
is simple. If $d=2$, then the linear system $|\L|$ has exactly 4 base points,
which are 4-torsion points in case $\L$ is symmetric. Furthermore, in this case
$\psi_\L$ is just a pencil of smooth (generically) irreducible  curves
of genus 3 [Ba]. If $d=3$, then $|\L|$ is base point free and $\psi_\L:
A\to\P_2$ is a 6-fold covering branched over a curve of degree 18
([LB], Example 10.1.5). If $d=4$, then $\psi_\L: A\to\P_3$
is generically birational on its image, a singular
octic surface ([BLvS], [LB] \S 10.5), and the geometry
of the situation is well understood. See [BLvS] and [LB],
Proposition 10.5.7 for
a complete description and for further details. Finally,
if $d\ge 5$, then Reider's theorem states that
$\psi_\L: A\to\P_{d-1}$ is an embedding if and only if there is no
elliptic curve $E$ on $A$ with $E\cdot c_1(\L)=2$. Furthermore,
the embedding is even projectively normal in case $d$ odd and
$d\ge 7$, or $d$ even and
$d\ge 14$ [La].\par

{\bf \S 1.3. Moduli.}

Let ${\A}_D$ denote the coarse moduli space of abelian varieties
of dimension $g$ with
polarization of type $D$. ${\A}_D$ is a $g(g+1)/2$-dimensional
quasi-projective variety. In case $g=2$ and $D=(1,d)$,
we'll write $d$ instead of $D$ in all our notation, so for instance
we'll write $\A_d$ instead of $\A_{(1,d)}$.
Similarly
${\A}^{lev}_D$ will denote
the moduli space ${\A}^{lev}_D$
of abelian varieties with a polarization of type $D$ and with canonical
level structure.

Relating these spaces, there is a
forgetful map $\pi_1:\A^{lev}_d\to \A_d$ which is a finite
dominant morphism, of degree $\sharp PSL(2,\boldz_d)=d(d^2-1)/2$
when $d$ is prime, and also a finite map induced by
dividing out the level structure
$\pi_2: \A^{lev}_d\to \A_1$, which is of degree
$d(d^4-1)/2$ when $d$ is prime (see [HKW], Proposition 1.21.).

Both $\A_D$ and $\A_D^{lev}$ can be described as certain quotients of the
Siegel upper half-space
$$\H_g=\{Z\in M_g({\bf C})\mid \hbox{$Z^t=Z$ and ${\rm Im}(Z)>0$}\}.$$
One can think of $\H_g$ as parametrizing abelian varieties with period matrix
$(Z,D)$, so that $\H_g$ is the moduli space of type $D$ abelian varieties
with a choice of a symplectic basis for the lattice,
namely the basis $\lambda_1,\ldots,\lambda_g,\mu_1,\ldots,\mu_g$ given
by the columns of the period matrix $(Z,D)$. Note that such data
also determines a canonical level structure, by using the
basis for $K(D)$ represented by
$\lambda_1/d_1,\ldots,\lambda_g/d_g,\mu_1/d_1,\ldots,\mu_g/d_g$.

There exists a universal family
$\X_D\rightarrow \H_g$ with $\X_D:=(\H_g\times {\bf C}^g)/\boldz^{2g}$,
where $\alpha\in \boldz^{2g}$ acts on $\H_g\times {\bf C}^g$ by
$$\alpha:(Z,v)\mapsto (Z, v+(Z,D)\alpha),$$
and there is a line bundle $\L$ on $\X_D$
such that $\L|_{X_Z}$ is a line bundle
of type $D$ and characteristic zero with respect to the decomposition
$Z\boldz^g\oplus D\boldz^g$, where $X_Z$ denotes the fibre of $\X_D\rightarrow
\H_g$ over $Z\in \H_g$. See [LB], \S 8.7 for details. Holomorphic
sections of $\L$ can be defined using the classical theta functions
$\vartheta\!\left [{c_1\atop c_2}\right ]$ for $c_1,c_2\in {\bf R}^g$. These
are holomorphic functions on $\H_g\times {\bf C}^g$. If $c_0,\ldots,
c_N$ are a set of representatives for the group $D^{-1}\boldz^g
/\boldz^g$, then the theta functions
$$\vartheta\!\left [{c_0\atop 0}\right ], \ldots,
\vartheta\!\left [{c_N\atop 0}\right ]$$
descend to give sections of $\L$ on $\X_D$, such that
their restrictions of $X_Z$ form a basis for
$H^0(\L|_{X_Z})$ for any $Z\in \H_g$.
Thus we can use these sections to
define a map $\X_D\rightarrow \P^N_{\H_g}$
in case $\L$ is generated by global sections.
The image of each abelian variety $X_Z$ under this map is Heisenberg
invariant, with the elements of $\H(D)$ acting by translation
on $X_Z$ by the corresponding elements of $K(D)$, as can be
seen from the relations
$$\vartheta\!\left [{c_{\nu}\atop 0}\right](Z,v+e_{\mu})=
\exp({2\pi i(c_{\nu}\cdot e_{\mu})})
\vartheta\!\left [{c_{\nu}\atop 0}\right](Z,v),$$ where
$e_{\mu}$ is the $\mu$-th standard basis vector, and
$$\vartheta\!\left [{c_{\nu}\atop 0}\right](Z,v+Ze_{\mu}/d_{\mu})=C_{\mu}
\vartheta\!\left [{c_{\nu}+e_{\mu}/d_{\mu}\atop 0}\right](Z,v),$$ where
$C_{\mu}$ is a constant which depends only on $\mu$.
The following proposition, which follows easily from
the above observations, demonstrates the significance of the moduli
spaces $\A_D^{lev}$:

\proclaim Proposition 1.3.1.  Let $X_Z$  and $X_{Z'}$ be two fibres of
$\X_D\rightarrow \H_g$, and suppose $\L|_{X_Z}$ and $\L|_{X_{Z'}}$
are both very ample. Then the images of $X_Z$ and $X_{Z'}$ coincide in
$\P^N$ if and only if there is an isomorphism between $X_Z$ and
$X_{Z'}$ preserving their canonical level structures.

Proof: We set $\L_Z:=\L|_{X_Z}$, $\L_{Z'}:=\L|_{X_{Z'}}$, and
let $\phi:X_Z\hookrightarrow \P^N$, $\phi':X_{Z'} \hookrightarrow \P^N$
be the maps induced by these bundles
using bases of classical theta-functions. Now, if
$\alpha:X_Z\rightarrow X_{Z'}$ is an isomorphism
of polarized abelian varieties, then there exists a
$c\in X_Z$ such that $\alpha^*\L_{Z'}=t_c^*\L_Z$. In particular there
is a linear automorphism $\beta:\P^N\rightarrow \P^N$ such that $\beta|_{X_Z}
=\alpha\circ t_{-c}$. On the other hand, $\alpha$ is an isomorphism of
level structures, say $b$ and $b^\prime$ respectively,
if and only if the composition $({b^\prime})^{-1}\circ \alpha \circ b$
$$K(D)\mapright{b} K(\L_Z) \mapright{\alpha} K(\L_{Z'})\mapleft{b'} K(D)$$
is the identity map. In other words, $\alpha$ is an isomorphism of
level structures if and only if
$$\alpha\circ t_{b(x)} \circ \alpha^{-1} = t_{b'(x)}\quad\hbox{for all
$x\in K(D)$,}$$
or equivalently that
$$\beta\circ t_{b(x)} \circ \beta^{-1} = t_{b'(x)}\quad\hbox{for all
$x\in K(D)$.}$$
Thus thinking of $\P^N=\P(W)$, and $\beta\in PGL(W)$, we deduce that
$\beta\in N(\H(D))/{\bf C}^*$, where $N(\H(D))$ is the normalizer of
$\H(D)$ in $SL^{\pm}(W)$. Now by [LB], Exercise 6.14,
there is an exact sequence
$$\exactm{K(D)}{N(\H(D))/{\bf C}^*}{Sp(D)},$$
where $Sp(D)$ is the group of symplectic automorphisms of $(K(D),e^D)$.
It follows that $\alpha$ is an isomorphism of level structures if and only
if $\beta$ induces the identity in $Sp(D)$, that is if and only if
$\beta\in K(D)\subseteq PGL(W)$, which is the case if and only if
$\phi'(X_{Z'})=\beta(\phi(X_Z))=\phi(X_Z)$.
$\bullet$

Thus, if $\L|_{X_Z}$ is very ample for general $Z\in \H_g$, there is
an open set $U\subseteq \A_D^{lev}$ over which we have a family
$\X_U\rightarrow U$, $\X_U\subseteq \P^N\times U$ of
embedded abelian varieties.
We will often make use of this family in the sequel. We should note, however,
that in general this is not the universal family over $\A_D^{lev}$.
Indeed, let $\psi:\H_g\rightarrow \A_D^{lev}$ be the quotient map
defining $\A_D^{lev}$ as a quotient of $\H_g$, and let $U'=\psi^{-1}(U)$.
Then the projective family $\X_U$ described above over $U$ pulls back
to the family $\X_D\rightarrow \H_g$ restricted to $U'$. However
$\X_U$ is a quotient of $\X_D$ which may identify two abelian
varieties over two different points of $\H_g$ via morphisms
which are not homomorphisms of abelian varieties, i.e., the
zero section of $\X_D\rightarrow\H_g$ is not preserved under taking
this quotient. Therefore the family $\X_U\rightarrow U$ is in fact a twist of
the universal family over $U$. The simplest case where this happens is
for $g=1$ and $D$ even (see [BaH] for further details).
\par\bigskip
\vfill\eject

\noindent
{\hd \S 2. Equations for Elliptic Curves and Abelian Surfaces.}

Let $A$ be an abelian variety of dimension $g$, $\L$ a symmetric line
bundle on $A$ of type $D=(d_1,\ldots,d_g)$. Choose a decomposition
for $\L^2$, with $K(\L^2)=K_1\oplus K_2$, with
$K_1\cong K_2\cong \boldz^g/2D\boldz^g$. This induces a decomposition on
$2K(\L^2)=K(\L)=2K_1\oplus 2K_2$. There is a canonical choice of basis
of $H^0(\L^2)$ given by $\{\vt_i^{\L^2}|i\in K_1\}$, and a canonical
choice of basis of $H^0(\L)$ given by $\{\vt_i^{\L}|i\in 2K_1\}$.

If $p_1$ and $p_2$ are the first and second projections of $A\times A$
onto $A$, and $\L$ and $\M$ are line bundles on $A$, then we write
$\L\boxtimes\M$ for $p_1^*\L\otimes p_2^*\M$. Similarly, if $s\in
\Gamma(\L)$ and $t\in\Gamma(\M)$ are sections, we write $s\boxtimes t$ for
$p_1^* s\otimes p_2^* t$.

\proclaim Theorem 2.1. The $(\#2K_1)\times(\#2K_1)$-matrix
$$M_{\L}=(\sum_{z\in Z_2}\vt^{\L^2}_{i+j+z}\boxtimes\vt^{\L^2}_{i-j+z})_{
2i,2j\in 2K_1}$$ has rank at most one on $A\times A$.
This notation means that $i$ (and ditto for $j$) runs through a subset $R$
of $K_1$ of representatives for $K_1/Z_2$.
Also  $Z_2=A_2\cap K_1\cong(\boldz/2\boldz)^g$, where
$A_2$ is the group of two-torsion points on $A$.

Proof: Define $\alpha:A\times A\rightarrow A\times A$ by $\alpha(a_1,a_2)
:=(a_1+a_2,a_1-a_2)$. $\alpha$ is an isogeny with $\ker(\alpha)
=\{(a,a)\mid a\in A_2\}$. By [LB], Lemma 7.1.1,
$$\alpha^*(\L\boxtimes\L)=\L^2\boxtimes\L^2.$$
Also,
$$K(\L^2\boxtimes\L^2)=(K_1\times K_1)\oplus(K_2\times K_2)$$
and
$$K(\L\boxtimes\L)=(2K_1\times 2K_1)\oplus(2K_2\times 2K_2),$$
and these decompositions are compatible with $\alpha$ in the sense
that $\alpha(K_i\times K_i)\cap K(\L\boxtimes \L)=2K_i\times 2K_i,$
by [LB], Lemma 7.1.2. Thus we are in a position to apply the Isogeny
theorem,
[LB], Theorem 6.5.1.
The K\"unneth isomorphism $H^0(\L\boxtimes\L) \cong H^0(\L)\boxtimes
H^0(\L)$ yields a basis of canonical theta-functions
$$\{\vt_{i'}^{\L}\boxtimes \vt_{j'}^{\L}\mid i',j'\in 2K_1\}=
\{\vt_{2i}^{\L}\boxtimes \vt_{2j}^{\L}\mid 2i,2j\in 2K_1\}.$$
Thus by the Isogeny theorem,
$$\alpha^*(\vt_{2i}^{\L}\boxtimes \vt_{2j}^{\L})
=\sum_{(a_1,a_2)\in\atop\alpha^{-1}(2i,2j)\cap (K_1\times K_1)}
\vt_{a_1}^{\L^2}\boxtimes\vt_{a_2}^{\L^2}.$$
If $\alpha(a_1,a_2)=(2i,2j)$ for $(a_1,a_2)\in K_1\times K_1$, then
we obtain
$$\alpha^*(\vt_{2i}^{\L}\boxtimes \vt_{2j}^{\L})
=\sum_{z\in Z_2}
\vt_{a_1+z}^{\L^2}\boxtimes\vt_{a_2+z}^{\L^2}.$$
Since the matrix $M'=(\vt_{2i}^{\L}\boxtimes\vt_{2j}^{\L})_{2i,2j\in 2K_1}$
obviously has rank one on $A\times A$, so does $\M_{\L}=\alpha^*M'$.
$\bullet$

\proclaim Corollary 2.2. If $A\cong E$ is an elliptic curve, $\L$ a symmetric
line bundle of degree $d\ge 2$ on $E$, and $E\subseteq\P^{2d-1}$ the
embedding given by the sections $\vt^{\L^2}_i$, $i\in K_1(\L^2)$, then
the $2\times 2$ minors of the matrix
$$M_d=(x_{i+j}y_{i-j}+x_{i+j+d}y_{i-j+d})_{2i,2j\in 2\boldz/2d\boldz}$$
vanish along $E\times E\subseteq\P^{2d-1}\times\P^{2d-1}$, where $x_i, y_j$,
with
$i,j\in\boldz/2d\boldz$, are coordinates for the first and second
$\P^{2d-1}$, respectively.\par

\noindent
{\it Example 2.3.} Consider an elliptic curve $E$ with $\L$ a
symmetric line bundle of degree 2. Embedding $E\times E$ into
$\Pthree\times\Pthree$ via $\L^2\boxtimes\L^2$, we find that the matrix
$$M_2=\pmatrix{ x_0y_0+x_2y_2&x_1y_3+x_3y_1\cr
x_1y_1+x_3y_3&x_2y_0+x_0y_2\cr}$$
has rank one on $E\times E$. Thus, if we substitute the coordinates of a point
$(y_0:y_1:y_2:y_3)\in E\subset\Pthree$, we obtain an equation for $E$
in $x_0,\ldots,x_3$
given by
$$\det M_2=0.$$
Applying $\sigma$ to this equation yields another equation, and these
two equations cut out $E$:
$$\eqalign{
(x_0y_0+x_2y_2)(x_2y_0+x_0y_2)-(x_1y_1+x_3y_3)(x_1y_3+x_3y_1)&=0\cr
(x_3y_0+x_1y_2)(x_1y_0+x_3y_2)-(x_0y_1+x_2y_3)(x_0y_3+x_2y_1)&=0.\cr}$$
Furthermore, since $(y_0:y_1:y_2:y_3)\in E$ and hence is a solution of
the above equations, we deduce that the closure of the union
in $\Pthree$ of all Heisenberg invariant elliptic normal curves
coincides with the surface $F'$ whose equation is
the determinant of the matrix obtained from $M_{2}$ by substituting
$x$'s for the $y$'s:
$$F':=\{x_{0}x_{2}(x_{0}^{2}+x_{2}^{2})-x_{1}x_{3}(x_{1}^{2}+x_{3}^{2})=0\}
\subset\Pthree.$$
It is easily seen that $F'$ is in fact projectively equivalent with the Fermat
quartic $F:=\{\sum_{i=0}^{3}x_{i}^{4}=0\}\subset\Pthree$. Finally, as
noted at the end of \S 1.3, we remark here that $F'$ is not birational
to the Shioda surface $S(4)$ (see [BaH] for a proof of this claim).

\proclaim Corollary 2.4. If $A$ is an abelian surface, and $\L$ a
symmetric line bundle
of type $(1,d)$, then $K_1=K_1(\L^2)\cong \boldz/2\boldz\oplus
\boldz/2d\boldz$,
$Z_2\cong \boldz/2\boldz\oplus\boldz/2\boldz$, and the matrix $\M_{\L}$
$$(
\vt_{i+j}^{\L^2}\boxtimes\vt_{i-j}^{\L^2}+
\vt_{i+j+(1,0)}^{\L^2}\boxtimes\vt_{i-j+(1,0)}^{\L^2}+
\vt_{i+j+(0,d)}^{\L^2}\boxtimes\vt_{i-j+(0,d)}^{\L^2}+
\vt_{i+j+(1,d)}^{\L^2}\boxtimes\vt_{i-j+(1,d)}^{\L^2})_{2i,2j\in 2K_1}$$
has rank at most one on $A\times A$.\par

\noindent
{\it Remark 2.5.} Note that for any abelian variety $A$ and choice of $y\in A$,
the $2\times 2$ minors of $M_{\L}$ yield quadratic theta relations, and one
might want to know if these theta relations are equivalent to the Riemann
theta relations. This is indeed the case, if one also considers all
translates of
these relations by the Heisenberg group (as in Example 2.3, where
we had to apply $\sigma$ to $\det M_2$ in order to obtain defining
equations for the elliptic curve $E$).

To see this, first recall that if $\rho\in \hat Z_2=\hom(Z_2,{\bf C}^*)$ is
a character of $Z_2$, $y\in K_1(\L^4)$, then we may define
$$\vt^{\L^4}_{y,\rho}:=\sum_{z\in Z_2} \rho(z)\vt^{\L^4}_{y+z}.$$
Then the Riemann theta relations (compare [LB] Theorem 7.5.2)
take the form, for $y,y_1,y_2\in K_1(\L^4)$,
with $y\equiv y_1\equiv y_2\mod 2K_1(\L^4)$ and $\rho\in \hat Z_2$,
$$\vt^{\L^4}_{y_1,\rho}(0)\sum_{z\in Z_2} \rho(z)\vt^{\L^2}_{y+y_2+z}
\otimes\vt^{\L^2}_{y-y_2+z}=
\vt^{\L^4}_{y_2,\rho}(0)\sum_{z\in Z_2} \rho(z)\vt^{\L^2}_{y+y_1+z}
\otimes\vt^{\L^2}_{y-y_1+z}.$$
See also [Mu3] Chapter II, Theorem 6.1,
for a slightly different version of Riemann's
theta identities  and [Mu3], p.223 ff, for a discussion of the
relationships between the various settings.

To obtain the above theta relations using the $2\times 2$ minors
of the matrix $M_{\L}$, we may proceed as follows.
Let $R\subseteq K_1(\L^2)$ be the subset in Theorem 2.1 being used to
represent $K_1(\L^2)/Z_2$. Then it is possible to find $i,j,i',j'\in R$
such that
$$\eqalign{ i+j \equiv y+y_2 &\mod Z_2,\cr
i'+j' \equiv y-y_2 & \mod Z_2,\cr
i+j' \equiv y+y_1 & \mod Z_2,\cr
i'+j \equiv y-y_1 & \mod Z_2.\cr }$$
Thus there exists $z_1,z_2,z_3,z_4\in Z_2$ such that
$$\eqalign{
i+j&= y+y_2+z_1,\cr
i'+j'&= y-y_2+z_2,\cr
i+j'&= y+y_1+z_3,\cr
i'+j&= y-y_1+z_4\cr}$$
hold in $K_1(\L^2)$.
Then the $2\times 2$ minor of $M_{\L}$ involving rows $i$ and $i'$, and
columns $j$ and $j'$ takes the form
$$\eqalign{m_{(i,j),(i',j')}
:=&
\left( \sum_{z\in Z_2} \vt^{\L^2}_{y+y_2+z}
\boxtimes \vt^{\L^2}_{i-j+z+z_1}\right)
\cdot
\left( \sum_{z\in Z_2} \vt^{\L^2}_{y-y_2+z}
\boxtimes \vt^{\L^2}_{i'-j'+z+z_2}\right)-\cr
&-
\left( \sum_{z\in Z_2} \vt^{\L^2}_{y+y_1+z}
\boxtimes \vt^{\L^2}_{i-j'+z+z_3}\right)
\cdot
\left( \sum_{z\in Z_2} \vt^{\L^2}_{y-y_1+z}
\boxtimes \vt^{\L^2}_{i'-j+z+z_4}\right).\cr}$$
Since $M_{\L}$ is rank $\le 2$, $m_{(i,j),(i',j')}$ is identically
zero on $A\times A$.

Suppose now $x\in K_2(\L^2)$. Then for $z,z'\in Z_2$,  by [LB], Proposition
6.4.2, there is a constant $C_x$ depending only on $x$ such that
$$\eqalign{t_x^*\vt^{\L^2}_{y+y_2+z}\cdot t_x^*\vt^{\L^2}_{y-y_2+z'}
&= C_x e^{\L^2}(y+y_2+z,x)\cdot e^{\L^2}(y-y_2+z',x)\cdot \vt^{\L^2}_{y+y_2+z}
\cdot \vt^{\L^2}_{y-y_2+z'}\cr
&= C_x e^{\L^2}(2y+z+z',x)\cdot
\vt^{\L^2}_{y+y_2+z}
\cdot \vt^{\L^2}_{y-y_2+z'}.\cr}$$
Now
$$\sum_{x\in K_2(\L^2)} e^{\L^2}(z+z',x) = \cases{0&$z+z'\not=0$,\cr
\# K_2(\L^2) &$z+z'=0$,\cr}$$
so
$${1\over \# K_2(\L^2)}\sum_{x\in K_2(\L^2)}
{1\over C_x e^{\L^2}(2y,x)}\cdot t_x^* \vt^{\L^2}_{y+y_2+z}\cdot
t_x^*\vt^{\L^2}_{y-y_2+z'}
=\cases{\vt^{\L^2}_{y+y_2+z}\cdot \vt^{\L^2}_{y-y_2+z}&$z=z'$,\cr
0&$z\not=z'$.\cr}$$
The same holds if $y_2$ is replaced with $y_1$. Thus, if we denote translation
by $(x,y)\in A\times A$ by $t_{(x,y)}$, then
$$\eqalign{n_{(i,j),(i',j')}
:=&
{1\over \# K_2(\L^2)}\sum_{x\in K_2(\L^2)}
{1\over C_x e^{\L^2}(2y,x)}\cdot t_{(x,0)}^* m_{(i,j),(i',j')}=\cr
=&
\sum_{z\in Z_2} (\vt^{\L^2}_{y+y_2+z}\cdot \vt^{\L^2}_{y-y_2+z}
\boxtimes \vt^{\L^2}_{i-j+z+z_1}\cdot \vt^{\L^2}_{i'-j'+z+z_2})\cr
&-
\sum_{z\in Z_2} (\vt^{\L^2}_{y+y_1+z}\cdot \vt^{\L^2}_{y-y_1+z}
\boxtimes \vt^{\L^2}_{i-j'+z+z_3}\cdot \vt^{\L^2}_{i'-j+z+z_4})\cr}$$
Next, using the fact that if $z\in Z_2\subseteq K_1(\L^2)$, then
$t_z^*\vt^{\L^2}_x=D_z\vt_{x+z}^{\L^2}$, where $D_z$ is
a constant depending only on $z$ (cf. [LB], Proposition 6.4.2), we see that
$$\eqalign{r_{(i,j),(i',j')}=&
\sum_{z'\in Z_2} {1\over{D_{z^\prime}}}
\rho(z') t^*_{(0,z')} n_{(i,j),(i',j')}\cr
=&
\left(\sum_{z\in Z_2} \rho(z)\vt^{\L^2}_{y+y_2+z}\cdot
\vt^{\L^2}_{y-y_2+z}\right)
\boxtimes\left(\sum_{z\in Z_2} \rho(z)\vt^{\L^2}_{i-j+z+z_1}\cdot
\vt^{\L^2}_{i'-j'+z+z_2}\right)-\cr
&-
\left(\sum_{z\in Z_2} \rho(z)\vt^{\L^2}_{y+y_1+z}
\cdot \vt^{\L^2}_{y-y_1+z}\right)
\boxtimes\left(\sum_{z\in Z_2} \rho(z)\vt^{\L^2}_{i-j'+z+z_3}
\cdot \vt^{\L^2}_{i'-j+z+z_4}\right)\cr
=&
\left(\sum_{z\in Z_2} \rho(z)\vt^{\L^2}_{y+y_2+z}
\cdot \vt^{\L^2}_{y-y_2+z}\right)
\boxtimes\left( \vt^{\L^4}_{y_1,\rho}(0)
\cdot \vt^{\L^4}_{y',\rho}\right)-\cr
&-
\left(\sum_{z\in Z_2} \rho(z)\vt^{\L^2}_{y+y_1+z}
\cdot \vt^{\L^2}_{y-y_1+z}\right)
\boxtimes\left( \vt^{\L^4}_{y_2,\rho}(0)
\cdot \vt^{\L^4}_{y',\rho}\right),\cr}$$
the  latter equality by the multiplication formula ([LB], Theorem 7.1.4),
where $y'$ is such that
$$y'=i-j-y_1+z_1=i'-j'+y_1+z_2=i-j'-y_2+z_3=i'-j+y_2+z_4.$$
The latter equalities are a consequence of our
choice for  $i,j,i'$ and $j'$. Finally, since
$\vt^{\L^4}_{y',\rho}$ is not everywhere zero, we can divide through
by this function and obtain the desired Riemann theta relation,
from the fact that $r_{(i,j),(i',j')}$ is identically zero on $A$.

\proclaim Lemma 2.6. Let $A$ be an abelian surface, and let $\M$ be
a line bundle of type $(1,2d)$ and characteristic zero with respect to some
decomposition. Then there exists an abelian surface $A'$ and a double
covering $f:A'\rightarrow A$ such that $f^*\M\cong\L^2$, where
$\L$ is a symmetric line bundle of type $(1,d)$ on $A'$. Furthermore,
the decomposition on $A$ for $\M$ induces a decomposition on $A'$
for $f^*\M$ such that $f^{-1}(K_1(\M))=K_1(f^*\M)$.

Proof:
An unbranched double cover of an abelian variety is determined by giving
a two-torsion element $\tau$ of $\Pic^0(A)=\hat A$.
We then obtain a commutative diagram
$$\matrix{
&&&&0&&0&&&&&\cr
&&&&\mapdown{}&&\mapdown{}&&&&\cr
&&&&K(f^*\M)&&K(\M)&&&&\cr
&&&&\mapdown{}&&\mapdown{}&&&&\cr
0&\mapright{}&G&\mapright{}&A'&\mapright{f}&A&\mapright{}&0&&\cr
&&&&\mapdown{\phi_{f^*\M}}&&\mapdown{\phi_\M}&&&&\cr
&&0&\mapleft{}&\hat A'&\mapleft{\hat f}&\hat A&\mapleft{}&\hat
G&\mapleft{}&0\cr
&&&&\mapdown{}&&\mapdown{}&&&&\cr
&&&&0&&0&&&&&\cr}$$
where $\hat G$ is the subgroup of $\hat A$ generated by $\tau$.
This shows that
$$\eqalign{K(f^*\M)&= \phi^{-1}_{f^*\M}(0)\cr
&= (\hat f\circ \phi_{\M}\circ f)^{-1}(0)\cr
&= (\phi_{\M}\circ f)^{-1}(\hat G).\cr}$$
Now, given that $\M$ is of type $(1,2d)$, we have an exact
sequence
$$\exact{
\boldz/2d\boldz\oplus\boldz/2d\boldz
}{\phi_{\M}^{-1}(\hat G)}
{\boldz/2\boldz}.$$
This sequence is split if $\phi_{\M}^{-1}(\tau)\cap A_2$ is non-empty.
Here $A_2$ denotes the set of two-torsion points of $A$.
If the sequence splits, we must have
$$f^{-1}(\phi_{\M}^{-1}(\hat G))=
(\boldz/2\boldz\oplus\boldz/2d\boldz)^{\oplus 2},$$
whence $f^*\M$ is of type $(2,2d)$.

Now we can write $A=V/\Lambda$ and $A'=V/\Lambda'$ for some $\Lambda'\subseteq
\Lambda$. Furthermore, $\hat A=V/\Lambda(\M)$. The decomposition
$\Lambda=\Lambda_1\oplus\Lambda_2$ induces a decomposition
$V=V_1\oplus V_2$. Set $\Lambda_i(\M)=\Lambda(M)\cap V_i$.

To construct the desired $A'$, choose a
$\tau\in {1\over 2}\Lambda_2$, $\tau\not\in\Lambda(\M)$, representing a
non-zero
two-torsion element of $\hat A$ in $\phi_{\M}(A_2)$. Then
$K_1(f^*\M)=\Lambda_1(\M)/(\Lambda'\cap \Lambda_1(M))=f^{-1}(K_1(\M))$.
In addition,
since $\M$ is of characteristic zero with respect to the decomposition,
so is $f^*\M$. Furthermore, $A_2'\subseteq K(f^*\M)$. Thus by
[LB], Exercise 6.12, $f^*\M=\L^2$ for some symmetric line bundle
$\L$ on $A'$.
$\bullet$

\proclaim Corollary 2.7. Let $A$ be an abelian surface and $\M$ be a line
bundle of type $(1,2d)$, $d\ge 2$, of characteristic zero with respect to
some decomposition on $A$. Let $\psi_{\M}:A\rightarrow \P^{2d-1}$ be the
map induced by the $\vt^{\M}_i$s, $i\in K_1(\M)$. Then the
$d\times d$-matrix
$$M_d=(x_{i+j}y_{i-j}+x_{i+j+d}y_{i-j+d})_{2i,2j\in 2\boldz/2d\boldz}$$
has rank at most two
on $(\psi\times\psi)(A\times A)\subseteq \P^{2d-1}\times
\P^{2d-1}$.

Proof: Let $f:A'\rightarrow A$ be the double cover in Lemma 2.6.
The decomposition on $A$ induces a decomposition on
$A'$ for $\L^2$ which is compatible with $f$.
By the Isogeny
theorem, thinking of $K_1(f^*\M)=\boldz/2\boldz\oplus K_1(\M)$,
$$f^*(\vt_a^{\M})=\vt_{(0,a)}^{\L^2}+\vt_{(1,a)}^{\L^2},$$
for $a\in K_1(\M)$. Thus
$$\eqalign{
(f\times f)^*&(\vt_{i+j}^{\M}\boxtimes \vt_{i-j}^{\M}+
\vt_{i+j+d}^{\M}\boxtimes \vt_{i-j+d}^{\M})_{2i,2j\in 2K_1(\M)}=\cr
=&\ \biggl(
(\vt^{\L^2}_{(0,i+j)}+\vt^{\L^2}_{(1,i+j)})\boxtimes
(\vt^{\L^2}_{(0,i-j)}+\vt^{\L^2}_{(1,i-j)})+\cr
&+
(\vt^{\L^2}_{(0,i+j+d)}+\vt^{\L^2}_{(1,i+j+d)})\boxtimes
(\vt^{\L^2}_{(0,i-j+d)}+\vt^{\L^2}_{(1,i-j+d)})\biggr)_{2i,2j\in 2K_1(\M)}\cr
=&\ M_{\L}+(\sigma_{(1,0)}\times 1_{A'})^*(M_{\L}),\cr}$$
where $\sigma_{(1,0)}:A'\rightarrow A'$ is translation by $(1,0)\in
K_1(\L^2)$, $f\times f: A'\times A'\to A\times A$ is the cartesian
product of $f$ with itself, and $\M_{\L}$ is the matrix of Corollary 2.4.
Since $M_{\L}$ has rank at most one, so does
$(\sigma_{(1,0)}\times 1_{A'})^*(M_{\L})$. Thus their sum has rank
at most two.
$\bullet$

\proclaim Corollary 2.8. Let $A$ be an abelian surface and let $\M$ be a line
bundle of type $(1,2d+1)$, $(d\ge 2)$, of characteristic zero with respect
to some decomposition on $A$. Let $\psi_{\M}:A\rightarrow
\P^{2d}$ be the map induced by the sections $\vt_i^{\M}$, $i\in K_1(\M)$. Then
the $(2d+1)\times (2d+1)$ matrix
$$M'_{d}=(x_{(d+1)(i+j)}y_{(d+1)(i-j)})_{i,j\in\boldz/(2d+1)\boldz}$$
has rank at most $4$ on
$(\psi\times\psi)(A\times A)\subseteq \P^{2d}\times
\P^{2d}$.

Proof:
If $f:A'\rightarrow A$ is an arbitrary double cover, then $\L=f^*\M$ must be
of type $(1, 4d+2)$ since $2d+1$ is odd. As in the proof of Lemma 2.6,
it is possible to choose the double cover and a decomposition
on $A'$ for $f^*\M$ which is compatible with $f$, in such a way so that
$f^{-1}(K_1(\M))=K_1(\L)=\boldz/(4d+2)\boldz$.
Therefore, by the Isogeny theorem,
$$f^*\vt_x^{\M} = \sum_{y\in K_1(\L)\atop f(y)=x} \vt_y^{\L}
=\vt^{\L}_{2x}+\vt_{2x+2d+1}^{\L},$$
where $x\in K_1(\M)\cong\boldz/(2d+1)\boldz$, and the indices in the
last sum are considered modulo $4d+2$.
Now
$$\eqalign{
&{(f\times f)}^*(\vt_{(d+1)(i+j)}^{\M}\boxtimes \vt_{(d+1)(i-j)}^{\M})=\cr
&=\
\vt^{\L}_{i+j}\boxtimes \vt^{\L}_{i-j}+
\vt^{\L}_{i+j+2d+1}\boxtimes \vt^{\L}_{i-j+2d+1}+
\vt^{\L}_{i+j+2d+1}\boxtimes \vt^{\L}_{i-j}+
\vt^{\L}_{i+j}\boxtimes \vt^{\L}_{i-j+2d+1},\cr}$$
since, thinking of $0\le i,j \le 2d$, we have
$4(d+1)(i+j)=2(i+j) \mod 4d+2$, whence $2(d+1)(i+j)=i+j$,
or $i+j+2d+1 \mod 4d+2$.
Thus
$$
{(f\times f)}^*(\vt_{(d+1)(i+j)}^{\M}\boxtimes \vt_{(d+1)(i-j)}^{\M})_{i,j
\in\boldz/(2d+1)\boldz}
=M_{\L}+(\sigma_{2d+1}\times 1_{A'})^* M_{\L}.$$
Since these two matrices have rank at most two, their sum has rank at most
four on $A'\times A'$, and hence also the original matrix has rank at
most four on $A\times A$.
$\bullet$

\proclaim Corollary 2.9. Let $E$ be an elliptic curve and let $\M$ be a line
bundle of degree $2d+1$, $d\ge 2$, of characteristic zero with respect
to some decomposition on $E$. Let $\psi_{\M}:E\rightarrow
\P^{2d}$ be the map induced by $\vt_i^{\M}$, $i\in K_1(\M)$. Then
the $(2d+1)\times (2d+1)$-matrix
$$M'_{d}=(x_{(d+1)(i+j)}y_{(d+1)(i-j)})_{i,j\in\boldz/(2d+1)\boldz}$$
has rank at most $2$ on
$(\psi\times\psi)(E\times E)\subseteq \P^{2d}\times
\P^{2d}$.

Proof: The proof is exactly the same as the proof of Corollary 2.8.
$\bullet$
\par\bigskip

We remark that similar results hold for abelian varieties of
dimension $g$ and polarization of type $(1,1,\ldots,1,n)$, but
now the above matrices will have much higher rank (depending exponentially
on $g$) on these varieties, and it is not clear if much information is
obtained.
\par\bigskip

\noindent
{\it Example 2.10.} Let $E$ be an elliptic curve endowed with a
symmetric line bundle $\L$ of degree 7. Embedding $E\times E$ into
$\P^{6}\times\P^{6}$ via the canonical basis of $\L\boxtimes\L$,
we find that the matrix $M'_{3}$ in Corollary 2.9 has rank
at most $2$ on $E\times E$. Therefore substituting in the matrix
the coordinates of a point $(y_{0}:y_{1}:\ldots:y_{6})\in E$, we
obtain a $7\times 7$-matrix $M'_{3}(y)$ with linear entries in
$x_{0}, \ldots,x_{6}$ whose $3\times 3$-minors vanish on
the elliptic curve normal curve $E\subset\P^{6}$. On the other side, if we
substitute for $(y_{0}:y_{1}:\ldots:y_{6})$ the coordinates of the
origin of $E$ , which is the only point of the intersection
$E\cap(\Ptwo)^-$, then the matrix $M'_{3}(y)$ becomes skew-symmetric.
Therefore, the $4\times 4$-pfaffians of $M'_{3}(y)$ provide then
35 quadratic polynomials vanishing on $E$,
which in fact generate the homogeneous ideal
of $E$. To see this recall first that $h^0(\I_E(2))=14$
since $E\subseteq\P^6$  is projectively normal, and that $H^0(\I_E(2))$
decomposes as an $\HHH_7$-module as the direct sum of two mutually
isomorphic 7-dimensional representations (of weight 2) of $\HHH_7$.
Thus in order to prove our claim we need to determine the linear span of the
$\tau$-invariant pfaffians of the skew-symmetric matrix $M'_{3}(y)$.
There are exactly four such pfaffians, namely
$$
\matrix {
x_{0}^{2}y_{2}y_{3}+x_{1}x_{6}y_{3}^{2}-x_{3}x_{4}y_{1}^{2} &
x_{0}^{2}y_{1}y_{3}-x_{1}x_{6}y_{2}^{2}+x_{2}x_{5}y_{1}^{2} \cr
x_{0}^{2}y_{1}y_{2}-x_{3}x_{4}y_{2}^{2}+x_{2}x_{5}y_{3}^{2} &
x_{1}x_{6}y_{1}y_{2}-x_{3}x_{4}y_{1}y_{3}-x_{2}x_{5}y_{2}y_{3}. \cr}
$$
Their matrix of coefficients in terms of the $\tau$-invariant quadratic
monomials $x_{0}^{2}$, $x_{1}x_{6}$, $x_{2}x_{5}$, $x_{3}x_{4}$ is,
up to column operations, the skew symmetric matrix
$$
\pmatrix{
0&y_{1}y_{2}& -y_{2}y_{3}& -y_{1}y_{3}\cr
-y_{1}y_{2}&0& -y_{3}^{2}& y_{2}^{2}\cr
y_{2}y_{3}& y_{3}^{2}& 0& -y_{1}^{2}\cr
y_{1}y_{3}& -y_{2}^{2}& y_{1}^{2}& 0\cr}.
$$
Since the entries of this matrix do not have any common zeroes, we
deduce that the $4\times 4$-pfaffians of the matrix $M'_{3}(y)$
span precisely a 14-dimensional subspace of $H^{0}(\O_{\Psix}(2))$
if and only if the point $(0:y_{1}:y_{2}:y_{3}:-y_{3}:-y_{2}:-y_{1})$ lies
on the quartic plane curve $K$ defined by the pfaffian of the
coefficient matrix:
$$K:=\{(y_{1}:y_{2}:y_{3})\mid
y_{1}^{3}y_{2}-y_{2}^{3}y_{3}-y_{3}^{3}y_{1}=0\}\subset(\Ptwo)^{-}.$$
This smooth quartic, which is in fact the unique $PSL(2,\boldz_7)$-invariant
of degree $\le 4$ in $(\Ptwo)^{-}$ and is projectively equivalent to
the Klein quartic
$$K^{\prime}:=\{(y_{1}:y_{2}:y_{3})\in\Ptwo\mid
y_{1}^{3}y_{2}+y_{2}^{3}y_{3}+y_{3}^{3}y_{1}=0\},$$
is therefore the isomorphic image in $\P^{6}$ of the modular curve $X(7)$,
represented here by the $0$-section of $S(7)$, the Shioda surface of
level 7 (compare [Kl], [Ve]). In particular this proves our claim about
the homogeneous ideal of $E$. The equations above for an elliptic
normal curve in $\P^6$ (up to certain linear combinations) were first
found by Klein (cf. [Kl], [KlF], [Ve]). Moreover, the above
setting generalizes easily to all elliptic normal curves of odd degree,
giving rise to equations similar to those in [KlF], p. 245--246, [Kl2],
and [Ve], Proposition 5.8. We refer also to [GP] for a detailed
discussion of the geometry associated with $X(11)$.\par

We remark further that $Sec(E)$, the closure of the chordal variety to
$E$, is in fact contained in the locus $V$ defined by the
$6\times 6$-pfaffians of the matrix $M'_{3}(y)$,
for the same choice as above of $y$ as origin of $E$. Restricting to a
2-dimensional subspace (e.g. $\{x_{0}=x_{1}=x_{3}=x_{5}=0\}$) it is
easily seen that $V$ has the expected codimension 3, whence also degree
14 by Porteous' formula.  Since, on the other hand $Sec(E)$ is an
irreducible 3-fold of degree 14 (compare Proposition 5.1)
we deduce that the $6\times 6$-pfaffians of the matrix $M'_{3}(y)$
indeed cut out $Sec(E)$.

Finally  let now $(y_{0}:y_{1}:\ldots:y_{6})$ be a general point
on the elliptic normal curve $E\subset\P^6$.
In particular, we may assume for its coordinates that
$y_{i}\ne 0$, for all $i\in\boldz_7$. Then the determinant of the matrix
$M'_{3}(y)$ is a non-zero septic polynomial $F$ (e.g. since the coefficient
of $x_{0}^{7}$ in $F$ is $-\prod_{i=0}^{6}y_{i}\ne 0$) which
vanishes on the variety $Sec_{3}(E)\subset\P^{6}$ of trisecant planes to
the elliptic normal curve $E$. On the other hand $Sec_{3}(E)$ is an
irreducible fivefold, hence an irreducible hypersurface in $\P^{6}$,
whose equation being also $\HHH_{7}$-invariant must have degree
divisible by 7. It follows necessarily that $Sec_{3}(E)=\{F=0\}$.
\par\bigskip

In the remainder of this section, we relate the above
approach to the main result of [EKS], which in turn enables us
to give a more geometrical explanation for the above  matrices.
Let $X\subseteq\P^{n}=\P(V^{*})$ be a non-degenerate,
reduced, irreducible scheme and assume we can write $\L=\O_{X}(1)=
\L_{1}\otimes \L_{2}$ for suitable line bundles $\L_{1}$ and
$\L_{2}$ on $X$. Suppose also that $V_{i}\subset H^{0}(\L_{i})$,
$v_{i}=\dim V_{i}$, $i=1,2$,
are linear series such that the image of $V_{1}\otimes V_{2}$
through the multiplication map
$$\mu:H^{0}(\L_{1})\otimes H^{0}(\L_{2})\to H^{0}(\L)$$
is contained in $V\subset H^{0}(\L)$. Then the homogeneous
ideal $I_{X}$ of $X\subset\P(V^{*})$ contains the $2\times 2$-minors
of the 1-generic $v_{1}\times v_{2}$-matrix $M$ with entries linear
forms on $V$ corresponding to the multiplication morphism
$V_{1}\otimes V_{2}\to V$. As usual, 1-generic means
that all generalized entries of $M$ are non-zero.
In particular, if $X\subseteq\P^{n}$ is linearly
normal and $D$ is an effective divisor on $X$ which moves in
a linear system of (affine) dimension $v$, and whose
linear span in $\P^{n}$ has codimension $w$, then $I_{X}$
contains the $2\times 2$ minors of a $v\times w$-matrix $M_{|D|}$
with linear entries.

Conversely, if the homogeneous ideal of $X$ contains the
$2\times 2$-minors of the 1-generic, $v_{1}\times v_{2}$-matrix
$M$ of linear forms, associated  to a pairing
$V_{1}\otimes V_{2}\to V$, then  we can recover line
bundles and linear series as above through the following
correspondence:
$$
\eqalign{
\L_{1}=&\im(M:\O_{X}^{\oplus v_{1}}\to \O_{X}^{\oplus v_{2}}(1))\cr
\L_{2}=&\im(M^{t}:\O_{X}^{\oplus v_{2}}\to \O_{X}^{\oplus v_{1}}(1))\cr
V_{1}=&\im(M:H^{0}(\O_{X}^{\oplus v_{1}})\to H^{0}(\L_{1}))\cr
V_{2}=&\im(M^{t}:H^{0}(\O_{X}^{\oplus v_{2}})\to H^{0}(\L_{2}))\cr
}$$

The main result of [EKS] asserts that if $X$ is a reduced, irreducible
curve of genus $g$ and $\L_{i}$, $i=1,2$, are line bundles on $X$
of degrees at least $2g+1$, nonisomorphic in case both have
degree $2g+1$ and $g>0$, then the $2\times 2$-minors
of the matrix defined above generate the
homogeneous ideal $I_{X}$ of $X$ embedded
via $\L_{1}\otimes \L_{2}$.  In particular, this provides a
way to write down the equations of an elliptic normal curve as the
$2\times 2$-minors of a matrix with catalecticant blocks (see
[EKS], \S 3 (c)).

In this context, it is easy  to determine to which splitting
of the embedding of an elliptic normal curve do correspond the
determinantal equations in Corollary 2.2. It is precisely
the choice  of a level  structure which accounts for the nice symmetric
form of the matrices involved. We'll discuss only the following case:

Let $E\subseteq \P^{2d-1}$ be a Heisenberg invariant
normal curve of degree $2d$, with $d\ge 3$ and $d\equiv 1\mod 2$, and let
$y=(y_0:\ldots:y_{2d-1})\in E\subset\P^{2d-1}$ be a point on
the curve. Recall now from \S 1.1 that the generators $\sigma$ and
$\tau$ of $\HHH_{2d}$ act on $E$, via the Schr\"odinger
representation, by translation with genuine
$(2d)$-torsion elements, say $\rho$ and $\nu$.

Now let $D$ be the divisor
$D=\sum_{i=0}^{d-1}\tau^{2i}(y^{\prime})$ of degree $d$ on $E$,
where
$y^{\prime}:=\sigma^{d}\tau^{d}(y)=
(y_{d}:-y_{d+1}:\ldots:-y_{0}:y_{1}\ldots)\in E$. On the other
hand, since
$\sum_{i=0}^{d-1}\tau^{2i}(y^{\prime})\equiv d\cdot y^{\prime}+
\sum_{i=0}^{d-1}(2i)\cdot\nu\equiv d(d-1)\nu+d\cdot y^{\prime}$
and $d\cdot y^{\prime}\equiv d\cdot\sigma^{d}\tau^{d}(y)
\equiv d\cdot y+d^{2}\cdot\rho+d^{2}\cdot\nu
\equiv d\cdot y+ d(\rho+\nu)$
in the group law of the elliptic curve $E$, we deduce that
$$\L_{1}:=\O_{E}(D)=\O_{E}(d\cdot y+\gamma),$$
where $\gamma= d\cdot(\rho+\nu)$ is
a non-trivial 2-torsion point on $E$. Moreover, since in the group law
${2k}d\cdot\rho\equiv 0$ for all $k\in\NN$ ,
the set $\{D_{k}:=\sigma^{2k}(D)\mid k\in\{0,\ldots,d-1\}\}$
is a basis of sections for the linear series $H^{0}(\O_{E}(D))$.

The linear span $\Pi_0$ of the divisor $D_{0}:=D$ is a $d$-secant
$\P^{d-1}$ to $E$, whose equations are
$$\eqalign{
\Pi_{0}&=\{x_{0}y_{0}+x_{d}y_{d}=x_{1}y_{1}+x_{d+1}y_{d+1}=\ldots=
x_{d-1}y_{d-1}+x_{2d-1}y_{2d-1}=0\}\cr
&=\{x_{i}y_{i}+x_{i+d}y_{i+d}=0\mid {2i\in
2\boldz/2d\boldz}\}.}$$ Similarly,
the linear span of the divisor $D_{k}$,  $k\in\{0,\ldots,d-1\}$,
is the $d$-secant $\P^{d-1}$ to the curve $E\subset\P^{2d-1}$ defined by
$\Pi_{k}:=\sigma^{2k}(\Pi_{0})$. The chosen elliptic curve $E\subset\P^{2d-1}$
is embedded by the line bundle $\L:=\O_{E}((2d)o_{E})$,  therefore we can
take as a complementary line bundle
$\L_{2}:=\L\otimes \L_{1}^{-1}=\O_E(\iota(D))$,
where $\iota$ is the Heisenberg involution. This time, a basis of sections for
$H^0(\L_2)$ is the set $\{\iota(D_k)\mid k\in\{0,\ldots,d-1\}\}$. The
linear span of a divisor $\iota(D_k)$ is $\Pi^\prime_k:=\iota(\Pi_k)$,
the  $d$-secant $\P^{d-1}$ to $E$ defined by
$$\Pi^\prime_k=\{x_{k+i}y_{k-i}+x_{k+i+d}y_{k-i+d}=0,\ {2i\in
2\boldz/2d\boldz}\}.$$
Given the above chosen bases for $H^0(\L_i)$, $i=1,2$, the
$d\times d$-matrix of linear forms corresponding to the
multiplication map $H^{0}(\L_{1})\otimes H^{0}(\L_{2})\to H^{0}(\L)$
has as its $(i,j)$-entry the equation of the hyperplane
$\span(\Pi_j, \Pi_i^\prime)$, namely $x_{i+j}y_{i-j}+x_{i+j+d}y_{i-j+d}$.
In other words, we obtain the matrix $M_d$ in Corollary 2.2, whose
$2\times 2$-minors generate the homogeneous ideal of the elliptic
curve $E\subset\P^{2d-1}$.\par\medskip

\noindent
{\it Example 2.11.} Let $E\subset\P^5$ be a Heisenberg
invariant elliptic normal curve, and let
$y=(y_0:\ldots:y_{5})\in E\subset\P^{5}$ be a point on it.
We've seen above that $E$ is cut out by the $2\times 2$-minors
of the matrix
$$
M_3(y)=\pmatrix{x_0y_0+x_3y_3 & x_1y_5+x_4y_2 & x_2y_4+x_5y_1\cr
x_1y_1+x_4y_4 & x_2y_0+x_5y_3 & x_3y_5+x_0y_2\cr
x_2y_2+x_5y_5 & x_3y_1+x_0y_4 & x_4y_0+x_1y_3\cr}
$$
The two nets, say  $\F_1$ and $\F_2$, of trisecant planes to $E\subset\P^5$
defined by generalized  rows or columns respectively, both trace out
in $\P^5$ the cubic  hypersurface $V_y=\{\det M_3(y)=0\}$. In
particular, $V_y$ depends only  on the chosen decomposition
of the line bundle $\O_E(6o_E)$, i.e., only on the point $y$, and not on
our choices for bases.  Analyzing further this example, we have the
following results due to Veneroni and Room (see [Ro], 7.11, $(ii)$ and
9.22)

\proclaim Proposition 2.12 (Room).
${\{V_y\}}_{y\in E}$ is a linear pencil. Furthermore, $Sec(E)$ is the complete
intersection of any two members in this pencil. In particular,  any cubic
hypersurface having multiplicity two along $E$ is a member of the above
pencil, and hence is determinantal.

Proof. Each cubic $V_y$ is singular along $E$, and thus by B\'ezout
contains any secant line to the elliptic curve. On the other hand,
it is easily seen that a trisecant plane to $E$ which meets a
plane of the net $\F_1$ ($\F_2$) in a point outside $Sec(E)$
is necessarily a $\P^2$ in the net $\F_2$ ($\F_1$, respectively).
In particular, two trisecant planes contained in different $V_y$'s
either meet on $Sec(E)$, or they are disjoint. Therefore
any two different cubics $V_y$ and
$V_{y^\prime}$ intersect properly, and so the claim
follows since $\deg Sec(E)=9$ (see Proposition 5.1).
$\bullet$

\proclaim Remark 2.13.
\item{$i)$} The matrix $M_3(y)$ is symmetric (up to row and column operations)
if and only if $y$ is a 2-torsion point of $E$. In each of these four cases
the $2\times 2$-minors of the matrix $M_3(y)$ cut out a Veronese surface in
$\P^5$ containing $E$; the cubic hypersurface $V_y$  is
the secant variety of this Veronese surface.
\item{$ii)$} $\sigma(V_y)=V_{y+(2-{\rm torsion})}$ and
$\tau(V_y)=V_{y+(2-{\rm torsion})}$. In particular,
$V_y$ is invariant under the action of the subgroup
$\HHH'\subseteq\HHH_6$ generated by $\sigma^2$ and
$\tau^2$. The four cubic hypersurfaces
$V_{y}$, $\sigma(V_{y})$, $\tau(V_{y})$ and
$\sigma\tau(V_{y})$ span only a pencil, whose base locus
is $Sec(E)$.\par

Proof. For part $i)$ observe that
$M_3(y)$ is symmetric if and only if the line bundles
$\L_1$ and $\L_2$ coincide, that is exactly
when $y$ is a 2-torsion point on $E$. See [SR]
for the assertions concerning the Veronese surfaces.
Part $ii)$ is easy and left to the reader.
$\bullet$

\par\bigskip
\noindent
{\hd \S 3. Degenerations of abelian surfaces and elliptic curves.}

We will need the following facts about degenerations of
elliptic curves and abelian surfaces.

\proclaim Definition. Let $E\subseteq \P^{n-1}$ be a Heisenberg
invariant elliptic normal curve of degree $n$. For a point $\tau
\in E$, $\tau$ not a 2-torsion point of $E$, the surface
$$S_{E,\tau}:=\bigcup_{P\in E} \langle P, P+\tau \rangle$$ is
called a {\it translation scroll},
where $\langle P,P+\tau\rangle$ denotes the line spanned by $P$ and
$P+\tau$.

\proclaim Theorem 3.1. Let $n\ge 5$.
Let $S_{E,\sigma}$ be a translation scroll,
with $\sigma\in E$ general.
Then there exists a flat family $\A\rightarrow\Delta$,
a point $0\in\Delta$, along with a Heisenberg invariant embedding
$\A\subseteq \P^{n-1}_{\Delta}$ such that $\A_0\cong S_{E,\sigma}$
and $\A_t$ is a non-singular abelian surface for $t\in\Delta$,
$t\not=0$.

Proof. This follows from the results of [DHS].
More precisely, let
$$\Omega_{\tau}=\pmatrix{
\tau_{11}&\tau_{12}\cr
\tau_{21}&\tau_{22}\cr
1&0\cr
0&n\cr}$$
be the normalized period matrix for an abelian surface
with a polarization of type $(1,n)$, where
$\tau=\tau^t$ is a $2\times 2$ symmetric complex matrix
with $\Im\tau >0$. Then [DHS], Proposition 14
and \S3.4, yield a family $\A_V$ over $V=\Spec {\bf C}[T_1,T_2]$
such that for general $(T_1,T_2)\in V$, the fibre is a smooth
abelian surface with $\tau=\pmatrix{\tau_{11}&\tau_{12}\cr
\tau_{12}&\tau_{22}\cr}$, for some $\tau_{11}$ and $\tau_{22}$,
and the fibre over $(0,T_2)$, $T_2\not=0$, is a translation
scroll associated to an elliptic curve
$$E={\bf C}/(\boldz n+\boldz\tau_{22})$$
for some $\tau_{22}$ depending on $T_2$, and each possible $\tau_{22}$
occurs. Furthermore the element in $E$ yielding the translation scroll
is given by $x=[\tau_{12}]\in E$. In this way one obtains
all possible translation scrolls. In [DHS], \S 4, for a given neighborhood
$\Delta$ of a point $(0,T_2)\in V$,
a line bundle $\L$ is constructed on $\A_{\Delta}$, along with
sections $\vt_0,\ldots,\vt_{n-1}\in H^0(\L)$ which are invariant
under the action of the Heisenberg group. It is then shown in \S 5
of [DHS] that these sections induce an embedding $\A_{\Delta}
\hookrightarrow \P^{n-1}_{\Delta}$ for general $\tau_{12}$, if
$n\ge 5$. This gives the desired result.
$\bullet$

We also need similar results about degenerations of elliptic curves
with level $n$ structure.

\proclaim Definition. We call
$$X(\Gamma_n)=\bigcup_{i\in\boldz_n} l_{i,i+1}=\bigcup_{i\in\boldz_n}
\langle e_{i}, e_{i+1}\rangle,$$
the {\it standard $n$-gon}
where $l_{i,i+1}$ is the line
where all coordinates are zero except for
the $i$-th and $(i+1)$-st in $\P^{n-1}$.
It is a cycle of $n$ lines.

\proclaim Theorem 3.2. There is a family $\E\subseteq\P^{n-1}_{\Delta}$
such that $\E_0=X(\Gamma_n)$ and $\E_t$ is a non-singular, $\HHH_n$-invariant
elliptic normal curve of degree $n$.

Proof: This is standard, accomplished by embedding a Tate curve
(see [DR], \S 1) in $\P^{n-1}$. See [Ste], Theorem 2.3, for the precise
equations of this degeneration.
$\bullet$

\proclaim Definition. We set $S_n\subseteq \P^{n-1}$ to be the closure
of the union of all $\HHH_n$-invariant elliptic normal curves in
$\P^{n-1}$.

Note that by \S 1.3, $S_n$ is an irreducible surface, and by Theorem 3.2
contains $X(\Gamma_n)$, the standard $n$-gon. If $n$ is odd, the normalization
of $S_n$ is the Shioda surface of level $n$.
\par\bigskip
\noindent
{\hd \S 4. Toric degenerations}

In this section we'll describe a class of toric
degenerations of abelian surfaces carrying a polarization of type
$(1,d)$, and also certain classes of combinatorially
defined ``varieties'' with trivial canonical sheaf. We start by
recalling some definitions.

As usual, a (lattice) polytope in ${\bf R}^{n}$ is
the convex hull of a finite subset of $\boldz^{n}$.

A finite (integral) polyhedral complex $\Delta$ is a finite set
of lattice polytopes in ${\bf R}^n$, such that any face of a polytope in
$\Delta$ is a polytope in $\Delta$, and such that any two of the
polytopes in $\Delta$ intersect in a face of each of them.
The polytopes in $\Delta$ will be called the faces of $\Delta$.
The maximal faces will be called facets. In particular, we
will denote by $F(\Delta)$ the set of facets of $\Delta$.
In the special case when  $\Delta$ consists of a single polytope $P$
we will denote by $\partial P$ the boundary complex, which is the
complex formed by all proper faces of $P$.

As a matter of notation, the $f$-vector of a $d$-dimensional
polyhedral complex $\Delta$ is the vector in $\NN^{d+1}$ given by
$f(\Delta)=(f_0,f_1,\ldots,f_d)$, where $f_k=f_k(\Delta)$
denotes the number of $k$-dimensional faces in $\Delta$. Finally
by the $f$-vector of a polytope $P$ we mean the
$f$-vector of its boundary complex $\partial P$.

Recall that if $\Delta$ is a simplicial complex on the vertex set $V=
\{v_0,\ldots v_n\}$, the corresponding {\it face (ring) variety}
$X(\Delta)$ in the sense of Stanley, Hochster and Reisner
(see [Ho] for more details) is the variety defined by the
ideal $I_{\Delta}$ in
the projective space $\P^n=\Proj {\bf C}[x_0,\ldots x_n]$, whose coordinates
correspond to the vertices of $\Delta$,
and the ideal $I_{\Delta}$ consists of monomials corresponding to the
simplexes (faces) not contained in  $\Delta$.
It is easily seen that the dimension of $X(\Delta)$ as a
projective variety is the
dimension of the simplicial complex $\Delta$.
The ideal $I_{\Delta}$ is generated by square free monomials,
and in fact the above construction describes a 1:1
inclusion-reversing correspondence between simplicial complexes
on the vertex set $V$
and ideals $I\subseteq (x_0,\ldots ,x_n)^2$ generated by square-free
monomials.

In the special case when the topological realization
$|\Delta|$ is a manifold, it follows
from the work of Reisner [Re], Hochster and Roberts [HR]
and Stanley [Sta] that
$X(\Delta)$ is a Gorenstein scheme. Moreover,
the following isomorphisms hold:
$$
H^i(X(\Delta),\O_{X(\Delta)}(n))=\cases{
H^i(\Delta,\CC)& for $n=0$ and all $i>0$;\cr
0& for all $n\ne 0$ and $0<i<{\rm dim}\, \Delta$,}$$
while
$$H^1(\P^n,\I_{X(\Delta)}(n))=\cases{
\widetilde H^0(\Delta,\CC)& if $n=0$;\cr
0& for all $n\ne 0$.\cr}$$
In particular, $X(\Delta)$ is projectively Cohen-Macaulay
if and only if $\Delta$ has no reduced cohomology below
${\rm dim}\,\Delta$, and $X(\Delta)$ is projectively
Gorenstein if and only if $|\Delta|$ is a homology sphere.
Moreover, it follows from [Ho] and [BE]
that the canonical bundle $\omega_{X(\Delta)}$
is 2-torsion, and moreover that it is trivial
if and only if $|\Delta|$ is an orientable manifold.

In terms of projective invariants $X(\Delta)$ has degree
$\deg X(\Delta)=f_{d}(\Delta)$ and arithmetic sectional genus
$\pi(X(\Delta))=df_d(\Delta)-f_{d-1}(\Delta)+1$, where
$d=\dim\Delta$.  Also, as expected,
$\chi(\O_{X(\Delta)})=\chi(\Delta)=\sum_{i=0}^d(-1)^if_i(\Delta)$.
In fact the Hilbert polynomial, which coincides with the Hilbert function for
strictly positive values, is completely determined by the
combinatorial data (see [Ho], or [BH] for more details).

A somewhat classical example of a face ring variety is
the following:

{\it Degenerations of elliptic normal curves.}
Let $\Gamma_{n+1}$ denote the triangulation of the circle $S^1$
corresponding to an $(n+1)$-gon, whose vertices are labeled in
a counterclockwise manner by $x_i$, $i\in\boldz_{n+1}=\boldz/(n+1)\boldz$.
Then $X(\Gamma_{n+1})\subseteq\P^n$ is a projectively
Gorenstein 1-dimensional scheme of degree $n+1$ and arithmetic
genus 1. In case $n=2$, $X(\Gamma_3)=\{x_0x_1x_2=0\}
\subseteq\P^2$ is just one of the triangles in the Hesse pencil.
To simplify notation in the general case we introduce the following:

\proclaim Definition. If $k,l\in\boldz_m$, we define the {\it distance}
between $k$ and $l$ to be
$$d(k,l)=\min\{|\bar k-\bar l|\mid  \hbox{$\bar k,\bar l\in\boldz$ representing
$k$ and $l$ respectively}\}.$$

It is now easy  to see that the homogeneous ideal
$I_{X(\Gamma_{n+1})}$, $n\ge 3$,  is generated by all quadratic
monomials $x_ix_j$,  where $i,j\in\boldz_{n+1}$ with
$d(i,j)\ge 2$. There are exactly ${(n^2-n-2)/ 2}$ such monomials.
It is also not hard to see that
$$X(\Gamma_{n+1})=\cup_{i\in\boldz_{n+1}} l_{i,i+1}\subseteq\P^n,$$
where $l_{i,i+1}=\langle e_i,e_{i+1}\rangle=
\{x_0=x_1=\ldots=x_{i-1}=x_{i+2}=\ldots=x_n=0\}$
is the line joining the vertices $e_i$ and $e_{i+1}$ of
the standard simplex. This fits with the notation of \S 3.
Finally we remark that $X(\Gamma_{n+1})$
is invariant under the action
of the extended Heisenberg group $\HHH^e_{n+1}$ via the
Schr\"odinger representation.

In the sequel, we will be interested in giving
a satisfactory description of
the ``secant varieties'' of the ``elliptic normal curves''
$X(\Gamma_{n+1})\subseteq\P^n$.\par\bigskip

{\it Cyclic polytopes.} These polytopes were first introduced by
Carath\'eodory and then rediscovered and studied by Gale, Motzkin and Klee
among others (see [Gr] for details and historical remarks).
The real algebraic curve $C_d\subseteq\RR^d$, parametrized by
$$t\mapsto x(t)=(t,t^2,\dots, t^d),\qquad t\in\RR,$$
is called the moment curve. The cyclic d-polytope $C(n,d)$
is defined as the convex hull  of $n$ distinct
points $x_i=x(t_i)$, $t_0<t_1<\ldots t_{n-1}$, $n\ge d+1$,
chosen on the moment curve $C_d$. Here are some relevant facts about
cyclic polytopes:

\item{-}  $C(n,d)$ is a $d$-dimensional simplicial polytope
(any $d+1$ points on $C_d$ being affinely independent). In
particular the boundary complex $\Delta(n,d):=\partial C(n,d)$ provides
a triangulation of the corresponding sphere:
$|\Delta(n,d)|\cong S^{d-1}$.

\item{-}  A $d$-uple $W$ of points in $V_{n,d}=\{x_0,\ldots,x_{n-1}\}$,
the set of vertices of $C(n,d)$, spans a facet of the cyclic polytope
if and only if any two vertices in $V_{n,d}\setminus W$ are
separated on the moment curve $C_d$ by an even number of points
of $W$. This criterion is called ``Gale's evenness condition''.
A similar description holds for the lower dimensional faces of
the cyclic polytope (see e.g. [BH], Theorem 5.2.11) and shows
that $C(n,d)$ is well defined as a combinatorial equivalence class
of polytopes.

\item{-} $C(n,d)$ is $[{d\over 2}]$-neighbourly, that is the
convex hull of any $j+1$ vertices of $C(n,d)$ is a face
of the polytope, for all $j\le [{d\over 2}]-1$. In particular,
$C(n,d)$ is a polytope without diagonals when $d\ge 4$.

\item{-} Gale's evenness condition above implies that
$$f_{n-1}(C(n,d))=\cases{
{n\over{n-k}}{{n-k}\choose k}& if $d=2k$;\cr
2{{n-k-1}\choose k}& if $d=2k+1.$\cr}$$
In particular, for a cyclic 4-polytope:

$$
f_{0}(C(n,4))=n,\quad f_{1}(C(n,4))={n\choose 2},
\quad f_{2}(C(n,4))=n(n-3),
\quad f_{3}(C(n,4))={1\over 2}n(n-3).$$

For $d=2$, the cyclic polytopes $C(n,2)$ are combinatorially
equivalent to the  $n$-gons $\Gamma_n$ introduced above.
When $d=4$, the boundary complexes of the
cyclic polytopes provide triangulations of the 3-sphere, and
therefore their associated face varieties will have the same
numerical type as a Calabi-Yau threefold. More precisely:

\proclaim Proposition 4.1.  $X(\partial C(n,4))\subseteq\P^{n-1}$,
$n\ge 5$, is a 3-dimensional projectively Gorenstein scheme
of degree $n(n-3)/2$ and arithmetic sectional genus $n(n-3)/2+1$,
with trivial Euler characteristic and trivial canonical sheaf.
Moreover\hfill\break
$i)$ $X(\partial C(n,4))$ is the ``secant variety''
of the ``elliptic normal'' curve $X(\Gamma_n)\subseteq\P^{n-1}$.
\hfill\break
$ii)$ $X(\partial C(5,4))=\{x_0x_1x_2x_3x_4=0\}$, while for $n\ge 6$
the homogeneous ideal $I_{X(\partial C(n,4))}$ is generated
by the cubic monomials
$$x_ix_jx_k,\quad\hbox{ with $d(i,j)\ge 2$, $d(i,k)\ge 2$ and $d(j,k)\ge 2$.}$$
There are $n(n^2-9n+20)/6$ such monomials.\par

Proof. The proof is easy and left to the reader. By ``secant variety'' to
$X(\Gamma_n)$ we mean the Zariski closure of
$$\bigcup_{x,y\in X(\Gamma_n)\atop x\not= y} \langle x, y \rangle,$$
where $\langle x,y \rangle$ is the line spanned by $x$ and $y$.
Statements $(i)$ and
$(ii)$ are immediate consequences of Gale's evenness condition
stated above.
$\bullet$

We remark here for later reference that the  coordinate ring
of $X(\partial C(n,2d))$, for $n\ge 6$ and $d\ge 2$, is an extremal graded
Gorenstein algebra in the sense of [Sch] (or a compressed or
extremally compressed Gorenstein
algebra of type $z^{d}$ in the sense of [FL]).  Let $A$ be a graded Gorenstein
$k$-algebra ($k$ a field) with $A_0\cong k$, and generated
by its degree one elements. Let $r=\dim_{\!k}(A_1)$ be its embedding
dimension, write $A$ as a quotient $k[x_1,\ldots,x_r]/I$
and define $d(I)$ (or $d(A)$) to be $\min\{ t\mid I_t\ne 0\}$.
Then the following inequality holds, see [Sch], [FL]:
$$ i(A)+ \dim A\ge 2 d(A)-1,$$
where  $i(A)$ is the index of regularity (that is
the degree from which
the Hilbert polynomial and the Hilbert function of $A$ start to agree),
and $\dim A$ is the Krull dimension. The Gorenstein algebra $A$ is called
extremal whenever equality holds in the above inequality.
In our case $i(I_{X(\partial C(n,2d))})=1$ while
$d(I_{X(\partial C(n,2d))})=d+1$.
Typical other examples of extremal graded
Gorenstein rings in (Krull) dimension 3 are the tangent cones of the simple
elliptic surface singularities.

We discuss in the sequel the possible topological types of the
face variety-like ``divisors'' on the 3-fold $X(\partial C(n,4))$.

\proclaim Proposition 4.2. If $\Delta$ is a triangulation
of a compact, connected, orientable 2-manifold $T_g$ which
can be realized as a simplicial subcomplex of the boundary
complex $\partial C(n,4)$ of the cyclic polytope $C(n,4)$,
then $g\le 1$. Moreover, in case $\Delta$ is a
triangulation of the torus $T_1$, then all vertices of the cyclic
polytope $C(n,4)$ are vertices of the triangulation, and the
1-skeleton of the triangulation has a Hamiltonian circuit
such that every triangle in $\Delta$ has exactly one edge
in common with this circuit.\par

Proof. Gale's evenness condition says that the convex
hull of four vertices of $C(n,4)$ is a facet of the polytope
if and only if they correspond to two pairs of neighbouring
points on the moment curve $C_4$. Similarly, three vertices
of $C(n,4)$ span a 2-face of the polytope if and only if
two of them are neighbours on the moment curve. For the first
part of the proposition observe that since there are at most
$m=f_0(\Delta)\le n$ pairs of vertices in $\Delta$ which are
adjacent on $C_4$ it follows that $f_2(\Delta)\le 2m$. Now $T_g$ is  a
manifold so each edge of $\Delta$ is contained in
exactly two triangles of $\Delta$. Therefore, since
$f_2(\Delta)=\chi(T_g)+f_1(\Delta)-m=2(1-g)+3f_2(\Delta)/2-m$,
we deduce that $f_2(\Delta)=2m-4(1-g)\le 2m$ which implies $g\le 1$.
We assume now that $g=1$, and therefore that equality holds in all the
above inequalities, that is $f_2(\Delta)=2m$ and $m=n=f_0(C(n,4))$
since the edges whose spanning vertices are adjacent on the
moment curve $C_4$ form a Hamiltonian circuit, in other words
a graph with $m$ vertices and $m$ edges which has no proper subgraphs.
But every edge of $\Delta$ is contained in exactly two
triangles, and thus each of the $2n$ triangles of $\Delta$ must
contain exactly one edge of the Hamiltonian cycle.
$\bullet$

{\it Remark 4.3.} In the case $g=1$, we can rephrase the statement
of the above proposition by saying that the face variety $X(\Delta)$ is
a ``translation scroll'' of the ``elliptic normal curve''
$X(\Gamma_{n})\subseteq X(\partial C(n,4))\subseteq\P^{n-1}$,
which corresponds to the Hamiltonian cycle.\par\smallskip

Let $\Delta$ be a triangulation of the torus $T_1$ and
denote as usual by $f_i=f_i(\Delta)$ the number of $i$-dimensional faces
of $\Delta$. The triplet $f(\Delta)=(f_0,f_1,f_2)$ will be called
the $f$-vector of the given triangulation. Each edge in the
triangulation is common to exactly two triangles, so $2f_1=3f_2$ and
hence $f_0=f_1-f_2+\chi(T_1)={1\over 3}f_1$. On the other hand
obviously $f_1\le {f_0\choose 2}$, so we deduce that $f_0\ge 7$.
Therefore a  triangulation of the torus $T_1$ has at least
7 vertices, and moreover the above formulae show that for such a
triangulation the graph of its 1-skeleton is necessarily $K_7$,
the complete graph on seven vertices. Such a triangulation
was first constructed in 1949 by Cs\'asz\'ar [Cs]. It is unique
up to isomorphism and has an automorphism group of order 42.
The dual graph of its 1-skeleton divides the torus in
the well known 7-colourable map (see [Wh] for more details).

Inspired by this construction we describe in the proof of the following
proposition a uniform series of triangulations for the torus $T_1$.

\proclaim Proposition 4.4. For each $n\ge 7$ there exists a
$\boldz_n\times\boldz_2$-invariant triangulation $\Delta_n$ of the
torus $T_1$, whose $f$-vector is $f(\Delta_n)=(n,3n,2n)$, and
which can be realized as a subcomplex in the cyclic polytope $C(n,4)$.
\par

Proof. Let $x_i=x(t_i)$, $i\in 0, 1,\ldots, n-1$, with $t_0<t_1<\ldots
t_{n-1}$,
denote the vertices of the cyclic polytope $C(n,4)$ on the moment curve $C_4$
and label them in a natural way by the elements of $\boldz_n$.
Assume first that $n\ge 8$. We define then $\Delta_n$ to be the
2-dimensional simplicial complex  whose faces are the triangles
$(x_i\,x_{i+1}\,x_{i+4})$ and $(x_i\,x_{i+3}\,x_{i+4})$, for $i\in\boldz_n$,
their edges and their vertices. These triangles are chosen such as to form
a tessellation of the usual representation of the torus as a rectangular
stripe with the opposite edges identified with an appropriate Dehn twist.
Namely we have the following diagram:
$$
\epsfbox{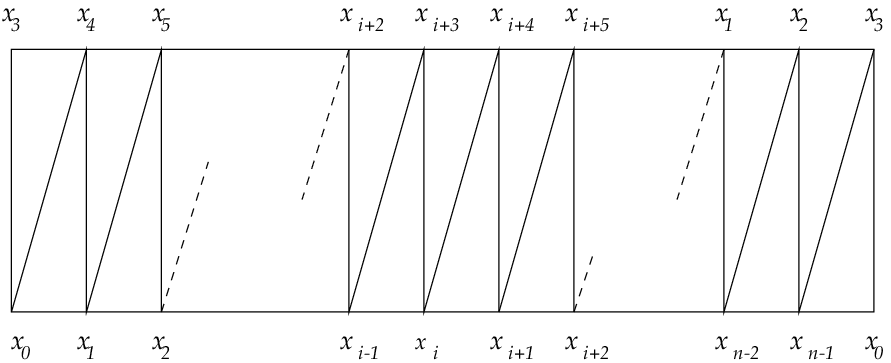}
$$
It is obvious now that $|\Delta_n|$ is homeomorphic to the torus $T_1$.
The $\boldz_n$ action on the triangulation is the one induced from
the natural action on the vertices. The extra $\boldz_2$ invariance is
under the transformation $i\mapsto -i\,(\mod n)$. A similar
discussion applies also in the case $n=7$, where $\Delta_7$ is taken
to be the 2-dimensional simplicial complex  whose faces are the triangles
$(x_i\,x_{i+1}\,x_{i+3})$ and $(x_i\,x_{i+2}\,x_{i+3})$, for $i\in\boldz_7$,
their edges and their vertices.
$\bullet$

\proclaim Corollary 4.5. The face variety $X(\Delta_n)\subseteq\P^{n-1}$ is
a locally Gorenstein, 2-dimensional  projectively normal
scheme of degree $2n$ and arithmetic sectional
genus $n+1$, with trivial canonical sheaf and irregularity $q=2$. In
particular,
$X(\Delta_n)$ has the same Hilbert polynomial as a smooth
$(1,n)$-polarized abelian surface. Moreover,
if $n\ge 13$, then  the homogeneous ideal $I_{X(\Delta_n)}$ is generated by
the quadratic monomials
$$x_ix_{i+2},\quad {\rm and}\quad x_ix_{i+5},\ldots, x_ix_{i+[{d\over 2}]},
\quad i\in\boldz_n,$$
where $[\dots]$ denotes the integral part. If $n=11$ or $12$, then
$I_{X(\Delta_n)}$ is generated by the previous quadrics
and the additional $n$ cubics
$$x_{i-4}x_ix_{i+4},\quad i\in\boldz_n,$$
while if $n=10$, $I_{X(\Delta_n)}$ is generated by the above quadrics
and the additional 10 cubics
$$x_{i-3}x_ix_{i+3},\quad i\in\boldz_{10}.$$
\par
Proof. In order to prove the statement about $I_{X(\Delta_n)}$,
$n\ge 13$, it is enough to
check that for each $(i,j,k)\in\boldz_n^3$, with $i\ne j$, $i\ne k$ and
$j\ne k$, either $(x_i\, x_j\, x_k)$ is a face of $\Delta_n$,
or one of the edges $(x_i\, x_j)$, $(x_i\, x_k)$, $(x_j\, x_k)$
is not an edge in $\Delta_n$. This is clear if either of $d(i,j)$,
$d(i,k)$, and $d(j,k)$ is two or at least 5. On the other hand,
if $n\ge 13$ and  all of $d(i,j)$, $d(i,k)$, and $d(j,k)$ are at
most 4, then it is easily seen that one of them is necessarily 2, if
$(x_i\, x_j\, x_k)$ is not a face in $\Delta_n$, while for $10\le n\le 12$,
the only non-faces $(x_i\, x_j\, x_k)$ with edges in $\Delta_n$ are the
ones listed.
$\bullet$

We now define
$$Q_i^{\lambda}:=\{x_i x_{i+2}+\lambda x_{i-1}x_{i+3}=x_j=0,\ {\rm for}\
j\in\boldz_n\setminus\{i,i+2,i-1,i+3\}\},$$ and then set
$$X_n^{\lambda}:= \cup_{i\in\boldz_n}Q_i^{\lambda},\quad \hbox{
for $\lambda\not=0$}.$$

Let $W_n$ be the one dimensional simplicial complex obtained by taking
the one-skeleton of $\Delta_n$ and removing the diagonals
$(x_i\, x_{i+4})$. We consider here $W_n$ as a simplicial subcomplex
of the boundary complex of the cyclic polytope $\partial C(n,4)$. Note
that the base locus $B$ of the family
$X_n^{\lambda}$ is in fact
$$B:=\bigcap_{\lambda\in k^*} X_n^{\lambda}= X(W_n).$$
This follows immediately from the fact that the ideal $I_{X(W_n)}$ is
generated by the monomials
$$\{x_ix_j\mid \hbox{$d(i,j)=2$ or $4$, or $d(i,j)\ge 5$}\}.$$
In particular, since the Hamiltonian cycle $\Gamma_n$ is a subcomplex of
$W_n$, we get $X(\Gamma_n)\subseteq X(W_n)$. The following was inspired and
is closely related to [ES], Proposition 4.8:

\proclaim Theorem 4.6. \item{$a)$}  The ideal
$I(X_n^{\lambda})$ is generated by the quadrics
$$\{x_ix_{i+2}+\lambda x_{i-1}x_{i+3} \mid i\in \boldz_n\}
\cup \{x_ix_j \mid d(i,j)\ge 5\}$$
for $\lambda \not=0$.
\item{$b)$} If $n=10$, then $x_{i-3}x_ix_{i+3}\in I(X_n^{\lambda})$
for all $i\in\boldz_{10}$, $\lambda\ne 0$,
while if $n=11$ or $12$, then $x_{i-4}x_i
x_{i+4}\in I(X_n^{\lambda})$ for all $i\in \boldz_n$.
\item{$c)$} The family $\{X_n^{\lambda}\mid \lambda\in k^*\}\subseteq
k^*\times \P^{n-1}$ extends uniquely to a flat family
$\{X_n^{\lambda}\mid \lambda\in {\bf A}^1\}\subseteq
{\bf A}^1\times \P^{n-1}$ over ${\bf A}^1$, and $X_n^0=X(\Delta_n)$. In
particular, $X_n^{\lambda}$ has the same Hilbert polynomial as a smooth
$(1,n)$-polarized abelian surface.

Proof: Let $J$ be the ideal generated by
$$\{x_ix_{i+2}+\lambda x_{i-1}x_{i+3} \mid i\in \boldz_n\}
\cup \{x_ix_j \mid d(i,j)\ge 5\}.$$
It is easy to see that $J\subseteq I(X_n^{\lambda})=
\bigcap_{i\in \boldz_n} I(Q_i^{\lambda})$.  For the purposes of
this proof, we define the width
of a monomial $x^I\in k[x_0,\ldots,x_{n-1}]$, where $I$ is a multindex, as
$$w(x^I):= \max\{d(i,j)\mid\hbox{$x_i|x^I$ and $x_j|x^I$}\}.$$
Suppose $f\in I(X_n^{\lambda})$. Since $I(X_n^{\lambda})\subseteq
I_{X(W_n)}$, it is clear that $f$ can have no terms of width 0 or 1,
and no terms of width three of the form $x_i^dx_{i+3}^e$. If $f$ has
a term of width two, say $x_i^{d_1}x_{i+1}^{d_2}x_{i+2}^{d_3}$, $d_1d_3
\neq 0$, then by using the relation $x_ix_{i+2}=-\lambda x_{i-1}x_{i+3}$,
we can replace this term, modulo $J$, with a term of width four or more. If $f$
has a term of width three, say $x_i^{d_1}x_{i+1}^{d_2}x_{i+2}^{d_3}
x_{i+3}^{d_4}$, $d_1d_4\neq 0$, then either $d_2$ or $d_3$ is non-zero
by the above discussion, and likewise we can replace this term
with a term of width
at least 4. Since $J$ contains all monomials of width at least 5, we thus find
that $f$ must be congruent modulo $J$ to a polynomial with terms of width
exactly 4.

If $n\ge 13$, then any width four term is of the form
$x_i^{d_1}x_{i+1}^{d_2}x_{i+2}^{d_3}
x_{i+3}^{d_4}x_{i+4}^{d_5}$, and if $d_3\neq 0$, we can again
replace this term with a term of width at least five, which thus is in
$J$. In case $d_3=0$, then restricting $f$ to the $\Pthree:=
\{x_j=0\mid j\in\boldz_n\setminus \{i,i+1,i+3,i+4\}\}$, $f$ must
vanish on the quadric $x_{i+1}x_{i+3}+\lambda x_ix_{i+4}$, and thus must
be a multiple of this quadric. Therefore we find that $f$ is congruent
modulo $J$ to a polynomial containing no terms of width less than or equal to
four and hence $f$ is zero. Thus $I(X_n^{\lambda})=J$ in this case.

If $n=10,11$ or $12$, the following minimal width four terms
could also appear:
$$\eqalign{ x_{i-3}x_i x_{i+3}&\quad\hbox{if $n=10$;}\cr
x_{i-4}x_i x_{i+4}&\quad\hbox{if $n=11,12$.}\cr}$$
Note that these are exactly the cubics appearing in Corollary 4.5.
But then
$$\lambda^2 x_{i-3}x_i x_{i+3}
=x_{i+3}\cdot (x_{i+1}x_{i+6})
-x_{i+6}\cdot (x_{i+1}x_{i+3}+\lambda x_i x_{i+4})
+\lambda x_{i}\cdot (x_{i+4}x_{i+6}+\lambda x_{i+3} x_{i-3}),$$
so $x_{i-3}x_ix_{i+3}\in J$.
Similarly, $x_{i-4}x_i x_{i+4}\in J$ for $i\in \boldz_n$, $n=11,12$,
via
$$\lambda x_{i-4}x_ix_{i+4}=x_{i+4}\cdot (x_{i-3}x_{i-1}+\lambda
x_{i-4}x_i) - x_{i-3}\cdot (x_{i-1}x_{i+4}).$$
This completes the proof of $a)$ and $b)$.
To prove $c)$, consider the family $\X_n$ in ${\bf A}^1\times\P^n$ defined
by the equations
$$\{x_ix_{i+2}+\lambda x_{i-1}x_{i+3} \mid i\in \boldz_n\}
\cup \{x_ix_j \mid d(i,j)\ge 5\}$$ along with the additional
cubics listed in part $b)$ if $10\le n \le 12$. Then the fibre
of this family over $\lambda\neq 0$ is $X_n^{\lambda}$, and by
Corollary 4.5, the fibre over $0\in {\bf A}^1$ is $X(\Delta_n)$. Thus
$\X_n$ is pure dimension three, and every component of $\X_n$ dominates
${\bf A}^1$. Therefore, by [Ha] \S III Proposition 9.7,
$\X_n\rightarrow {\bf A}^1$ is flat.
$\bullet$

\noindent
{\it Remark 4.7.} The smooth part $X(\Gamma_n)_{sm}$ of
$X(\Gamma_n)$ has a natural group structure, isomorphic to
$\Gm\times \boldz/n\boldz$. As such, it is still possible to
define a translation scroll for $X(\Gamma_n)$ by taking
a given $\tau\in X(\Gamma_n)_{sm}$, and letting $S_{X(\Gamma_n),\tau}$
be the Zariski closure of
$$\bigcup_{x\in X(\Gamma_n)_{sm}} \langle x, x+\tau \rangle.$$
Then it is not difficult to see that the $X_n^{\lambda}$
defined above are translation scrolls for $X(\Gamma_n)$,
and that the pencil $X_n^{\lambda}$ traces out $Sec(X(\Gamma_n))=X(\partial
C(n,4))$, with appropriate definitions for $\lambda=\infty$.

\par\bigskip
\noindent
{\hd \S 5. Secant varieties of elliptic normal curves.}

Let $E\subseteq\P^{n-1}$
be an elliptic normal curve of degree $n$ with a level $n$ structure.
This is equivalent to saying that the embedding of $E$ is
$\HHH_n\rtimes\langle\iota\rangle$-equivariant under the Schr\"odinger
representation,
where $\HHH_n$ acts on $E$ via translation with $n$-torsion points
while $\iota$ acts as the natural involution given by negation. We'll denote by
$Sec(E)$ the closure of the
chordal variety to $E$. We have

\proclaim Proposition 5.1. Let $E\subseteq\P^{n-1}$
be an elliptic normal curve of
degree $n$. Then
\item{$(i)$} $Sec(E)$ is an irreducible threefold of degree $n(n-3)/2$.
\item{$(ii)$} $Sec(E)$ is non-singular outside of $E$, and is singular
along $E$ with multiplicity $n-2$.
\item{$(iii)$} A natural desingularization of $Sec(E)$ is
given by $\pi:\P^1_{S^2E}\to Sec(E)$,
where
$$\P^1_{S^2E}=\{ (p,\{e_1,e_2\})\mid e_1,e_2\in E,
p\in\span(e_1,e_2)\}\subseteq
\P^{n-1}\times S^2E,
$$
and where $\pi$ is the projection on the first factor. $\P^1_{S^2E}$ has the
structure of a $\Pone$-bundle over the elliptic scroll $S^2E$.
\item{$(iv)$} $\omega_{Sec(E)}\cong \O_{Sec(E)}$, while
$h^1(\O_{Sec(E)})=h^2(\O_{Sec(E)})=0$.

Proof:
$i)$ $Sec(E)$ is clearly an irreducible threefold. To compute its degree,
take a general $L\cong\P^{n-4}\subseteq\P^{n-1}$ and project $E$ to
$\Ptwo$ from $L$. The number of nodes of the projection of $E$ will
be exactly ${1\over 2}(n-1)(n-2)-1=n(n-3)/2$, and this is precisely
the number of secants of $E$ which $L$ meets. This is also the degree of
$Sec(E)$.

$ii)$ (See [ADHPR1], [Seg] in the case $n=5$.)
To see that the multiplicity of $Sec(E)$ along $E$ is $n-2$, take a
linear space $L\cong \P^{n-4}$, $L\subseteq\P^{n-1}$, which contains
a point $p$ of $E$, so that $L$ is neither tangent to $E$ nor contains
a secant of $E$. Choosing $L$ generally, we can assume that $L$ meets
$Sec(E)$ transversally at a finite number of smooth points outside of $p$.
Now, two secants or tangents of $E$ cannot meet at a point outside of $E$.
If they do, then there is a plane $M$ containing both, and the linear
system of hyperplanes containing $M$ yields a linear system on $E$
residual to $M\cap E$ which is of dimension $n-4$ and degree $\le n-4$, a
contradiction for an elliptic curve. Thus there is a unique secant or tangent
line to $E$ through each point of $L\cap Sec(E)$ outside of $p$. Projecting
$E$ from $L$ yields a curve in $\Ptwo$ of degree $n-1$, which must then
have
${1\over 2}(n-2)(n-3)-1$ nodes. Thus $L$ intersects $Sec(E)$ in
${1\over 2}(n-2)(n-3)-1$ points outside of $p$, and so $p$ has multiplicity
$$\deg Sec(E)-\left[{1\over 2}(n-2)(n-3)-1\right]=n-2.$$
We will show that $Sec(E)$ is non-singular outside of $E$ after showing
$iii)$.

$iii)$ It is clear that $\Pone_{S^2E}$ is a $\Pone$-bundle over $S^2E$ via
the second projection $f:\Pone_{S^2E}\to {S^2E}$,
and thus is non-singular. Therefore $\pi$ is a
desingularization of $Sec(E)$.

To show that $Sec(E)$ is singular only
along $E$, it is enough to show that $\pi$ is an isomorphism outside of
$\pi^{-1}(E)$. Since no two bisecants of $E$ pass through the same
point of $\P^{n-1}$ outside of $E$, $\pi:\Pone_{S^2E}\setminus\pi^{-1}(E)
\rightarrow Sec(E)\setminus E$ is one-to-one, thus to check
that $\pi$ is an isomorphism,
we need to consider the maps on tangent spaces. Let
$x_1,x_2\in E$ be distinct points, and let $U_1,U_2$ be disjoint
small open neighborhoods of $x_1, x_2$, respectively, with local coordinates
$u_1$ and $u_2$. These yield a neighborhood of $\Pone_{S^2E}$ of the form
$U_1\times U_2\times {\bf A}^1$, with coordinates $(u_1,u_2,t)$. Let
${\bf A}^{n-1}\subseteq \P^{n-1}$ be an affine subspace containing $U_1$
and $U_2$, and let $\phi_1:U_1\rightarrow {\bf A}^{n-1}$,
$\phi_2:U_2\rightarrow {\bf A}^{n-1}$ be the embeddings of $U_1$
and $U_2$, respectively, in ${\bf A}^{n-1}$. Then we
may identify the restriction of $\pi$ with the map
$\psi:U_1\times U_2\times {\bf A}^1\rightarrow {\bf A}^{n-1}$
defined by
$$\psi(u_1,u_2,t)= t\phi_1(u_1)+(1-t)\phi_2(u_2).$$
For fixed $u_1$ and $u_2$, this yields the secant line joining $\phi_1(u_1)$
and $\phi_2(u_2)$. Now for a given point $(u_1,u_2,t)\in
U_1\times U_2\times{\bf A}^1$, the image $\psi_\ast(T)$ of the tangent space
$T=T_{U_1\times U_2\times {\bf A}^1,(u_1,u_2,t)}$ at that point is spanned by
the rows of the jacobian matrix
$$\pmatrix{\partial\psi/\partial u_1 (u_1,u_2,t)\cr
\partial\psi/\partial u_2 (u_1,u_2,t)\cr
\partial\psi/\partial t (u_1,u_2,t)\cr}
=\pmatrix{ t\phi_1'(u_1)\cr
(1-t)\phi_2'(u_2)\cr
\phi_1(u_1)-\phi_2(u_2)\cr}.$$
Therefore, if $t\not= 0,1 $ (values where $Sec (E)$ is singular),
then $\psi_*(T)$ contains
the tangent vectors to $E$ at $u_1$ and $u_2$, and hence $\psi_*(T)$
is tangent to $E$ at $u_1$ and $u_2$. If $\psi_*(T)$ had dimension less
than $3$, then this would yield a plane which intersects $E$ at least four
times counted with multiplicities, but again this is not possible.
Thus $\psi_*:T\rightarrow T_{\P^{n-1},\psi(u_1,u_2,t)}$ is
injective, and thus $Sec (E)$ is nonsingular away from $E$ at points
on secant lines which are not tangent lines.

To deal with points on tangent lines, we need suitable coordinates on
$S^2E$. Let $x\in E$, and let $U$ be an open neighborhood of $x$ with
coordinate $u$, with $u=0$ at $x$ for convenience. Let $\phi:U
\rightarrow {\bf A}^{n-1}$ be the embedding of $U$ as before. Consider
now $\psi:U\times U\times {\bf A}^1\rightarrow {\bf A}^{n-1}$ defined
by
$$\psi(u,v,t)=\cases{
{1\over 2}(\phi(u)+\phi(v))+t\left ( {\phi(v)-\phi(u)\over
v-u} \right )& $v\not= u$;\cr
\phi(u)+t\phi'(u)& $v=u$.\cr}$$
$\psi$ is clearly holomorphic and its image is contained in
the secant scroll. We can expand $\phi$ in  a Taylor series,
$$\phi(u)=\phi(0)+\phi'(0)u+{\phi''(0)\over 2} u^2+\cdots,$$
and then  write $\psi$ as
$$
\eqalign{
\psi(u,v,t)=&\ \phi(0)+{1\over 2}\left(\phi'(0)(u+v)+{\phi''(0)\over 2}
(u^2+v^2)+{\phi^{(3)}(0)\over 6} (u^3+v^3)+\cdots\right)+\cr
&+t\left ( \phi'(0)+{\phi''(0)\over 2} (u+v)+{\phi^{(3)}(0)\over 6}
(u^2+uv+v^2)+\cdots\right ).\cr
}$$
Now, on $S^2U$, we can use the symmetric coordinates $s_1=u+v$ and $s_2=uv$.
Since $\psi$ is symmetric with respect to interchanging $u$ and $v$,
$\psi$ descends to a function on $S^2U\times {\bf A}^1$, which we
may write as
$$
\eqalign{
\psi(s_1,s_2,t)=&\
\phi(0)+{1\over 2}\left(\phi'(0)s_1+{\phi''(0)\over 2}
(s_1^2-2s_2)+{\phi^{(3)}\over 6} (s_1^3-3s_1s_2)+\cdots\right)+\cr
&+t\left ( \phi'(0)+{\phi''(0)\over 2} s_1+{\phi^{(3)}(0)\over 6}
(s_1^2-s_2)+\cdots\right ).}
$$
As before, if $T$ is the tangent space of $S^2U\times {\bf A}^1$
at $(0,0,t)$, $t\not=0$, we can write $\psi_*(T)$ as the span of the
rows of the jacobian of $\psi$:
$$\pmatrix{\partial\psi/\partial s_1 (s_1,s_2,t)\cr
\partial\psi/\partial s_2 (s_1,s_2,t)\cr
\partial\psi/\partial t (s_1,s_2,t)\cr}
=\pmatrix{ {1\over 2}\phi'(0)+t{\phi''(0)\over 2}\cr
-{\phi''(0)\over 2}-t{\phi^{(3)}(0)\over 6}\cr
\phi'(0)\cr}.$$
So given that $t\not=0$, $\psi_*(T)$ is spanned by $\phi'(0)$,
$\phi''(0)$ and $\phi^{(3)}(0)$. Thus $\psi_*(T)$ intersects $E$ at
$x\in U$ with multiplicity four, and hence, as before, must be three
dimensional. Thus we have shown that $Sec(E)$ is smooth at all points in
$Sec (E)\setminus E$.

$iv)$ There is a natural inclusion $E\times E \subseteq \Pone_{S^2E}$,
given by $(e_1,e_2)\mapsto (e_1,\{e_1,e_2\})$. Thus the map
$\pi:\Pone_{S^2E}\rightarrow \P^{n-1}$ contracts
$E\times E\subseteq \Pone_{S^2E}$ to $E\subseteq \P^{n-1}$, via projection
onto the first coordinate. If $\widetilde X$ is the normalization of
$Sec(E)$, then
$\widetilde X$ is isomorphic precisely to the threefold obtained
by contracting $E\times E$ inside of $\Pone_{S^2E}$ in this manner.
Let $p:\Pone_{S^2E}\rightarrow\widetilde X$ be the contraction.
$\widetilde X$ then has a curve of simple elliptic singularities,
so $\widetilde X$ is Gorenstein and $K_{\Pone_{S^2E}}\sim p^* K_{\widetilde X}
-E\times E$. If $E\times E\sim -K_{\Pone_{S^2E}}+D$ for some divisor
$D$, adjunction tells us that
$$\eqalign{0=&\ K_{E\times E}\cr
=&\ (E\times E)\cdot (-K_{\Pone_{S^2E}}+D+K_{\Pone_{S^2E}})\cr
=&\ (E\times E)\cdot D.\cr}$$
It is easy to see that the restriction map $\Pic(\Pone_{S^2E})
\rightarrow \Pic(E\times E)$ is injective.
Thus $D\sim 0$ and $E\times E$ is an anticanonical divisor in the linear system
$|-K_{\Pone_{S^2E}}|$.  Thus
$$K_{\Pone_{S^2E}}\sim p^*K_{\widetilde X}+K_{\Pone_{S^2E}},$$
so $\K_{\widetilde X}\sim 0$. Finally using intersection theory in
$\Pone_{S^2E}$, one can see that $\widetilde X$ is singular with multiplicity
$n-2$  along $E$, and these singularities are locally of the form
$$\hbox{(simple elliptic surface singularity of multiplicity $n-2$)}
\times \hbox{curve}.$$
Thus the Zariski tangent space at each singular point of $\tilde X$ has
dimension $n-1$, and also the Zariski tangent space of each singular point
of $Sec (E)$ is of dimension $n-1$ (since such a tangent space contains the
cone over $E$ with vertex a point on $E$), and hence $\widetilde
X\cong Sec (E)$. Therefore $\omega_{Sec(E)}\cong \O_{Sec(E)}$,
as desired.

To conclude, by Serre duality it is enough to show that
$H^1(\O_{Sec(E)})=0$. Now $H^1(\O_{Sec(E)})$ is the Lie algebra
of $\Pic^0(Sec(E))$, which is reduced since we are in characteristic
zero. Thus it is enough to show that $\Pic(Sec(E))$ is a discrete group,
since then $H^1(\O_{Sec(E)})=0$. To show this, we use the map
$$\pi^*:\Pic(Sec(E))\rightarrow \Pic(\Pone_{S^2E}),$$
which is injective, and whose image is contained in the subgroup
$$\PP=\{\L\in\Pic(\Pone_{S^2E})\mid\hbox{$\L|_{E\times E}= \rho_1^*\M$ for
some $\M\in \Pic(E)$}\},$$
where $\rho_1:E\times E\rightarrow E$ denotes the first projection. So it is
enough to check that $\PP$ is discrete. We will do this by showing that
if $D,D'\in \PP$ are two divisors which are algebraically equivalent on
$\Pone_{S^2E}$, then $D\sim D'$. Indeed $D-D'\in \Pic^0(\Pone_{S^2E})$,
and any element of $\Pic^0(\Pone_{S^2E})$ can be written as a multiple of
$f^*(C_0-C_0')$, where $f:\Pone_{S^2E}\rightarrow S^2E$ is the second
projection, and
$C_0, C_0^\prime\subseteq S^2E$ are images of  fibres of $\rho_1$
in $S^2E$. But
$$f^*(C_0-C_0')|_{E\times E}= E\times \{p_1\}-E\times \{p_2\}
+\{p_1\}\times E-\{p_2\}\times E,$$
and this is not trivial on the fibers of the first projection $\rho_1:
E\times E\rightarrow E$ unless $p_1=p_2$, in which case $D\sim D'$.
Thus if $D-D'\in \PP$, then $D\sim D'$; and therefore
$\Pic(Sec(E))$ is discrete, and $H^1(\O_{Sec(E)})=0$.
$\bullet$

We now give the equations for these secant varieties. In this direction
we mention also Ravi's results [Rav] which in the special case of an elliptic
normal curve $E\subset\P^{n-1}$ say that if $M$ is the $m\times (n-m)$-matrix
with linear entries corresponding to a splitting of the embedding as
a tensor product of two line bundles of degrees $m$ and $n-m$, respectively,
with $4\le m\le n-4$, then the $3\times 3$-minors of $M$ set-theoretically
define $Sec(E)$.

\proclaim Theorem 5.2. Let $E\subseteq \P^{2d-1}$ be a Heisenberg invariant
normal curve of degree $2d$, $d\ge 3$. If $(y_0:\ldots:y_{2d-1})\in E$
is a general point, then the $3\times 3$ minors of the matrix
$$M_d=(x_{i+j}y_{i-j}+x_{i+j+d}y_{i-j+d})_{i,j\in 2\boldz/2d\boldz}$$
vanish on $Sec(E)$. In particular, if $d=3$, then $\det M_d$ and
$\sigma(\det M_d)$ cut out $Sec(E)$, and if $d\ge 4$, then for general $E$,
the ideal
of $3\times 3$ minors of $M_d$ cuts out $Sec(E)$.

Proof: By Corollary 2.2, the $2\times 2$ minors of $M_d$ vanish on $E$. It is
easy to see that the variety given by the determinant of a $3\times 3$ matrix
of forms is singular along the vanishing locus of the $2\times 2$ minors. Thus
each $3\times 3$ minor of $M_d$ determines a cubic hypersurface
which is singular along $E$. Any secant line to $E$ intersects this cubic
hypersurface in at least four points, counted with multiplicity, and so
is contained in this hypersurface. Thus each $3\times 3$ minor of $M_d$
vanishes along $Sec(E)$.

Now if $d=3$, it is easy to check that $\det M_3$ and $\sigma(\det M_3)$ are
independent cubics, and hence $\det M_3=
\sigma(\det M_3)=0$ determines a three dimensional complete intersection of
degree $9$. Since $Sec(E)$ is also degree 9, these two varieties must
coincide (see Proposition 2.12 and Remark 2.13, $ii)$).

Define now $\SS\subseteq\P^{2d-1}\times S_{2d}\subseteq \P^{2d-1}\times
\P^{2d-1}$, with $x$ coordinates on the first $\P^{2d-1}$ and
$y$ coordinates on the second $\P^{2d-1}$, as the variety defined
by the $3\times 3$ minors of the matrix $M_d$. (See \S 3 for
the definition of $S_{2d}$). If $E\subset\P^{2d-1}$ is a
Heisenberg invariant elliptic curve and $y\in E$, then
$y\in S_{2d}$ and we have already
observed that the fibre $\SS_t$ of the second projection $\SS\rightarrow
S_{2d}$ contains $Sec(E)$. On the other hand, if there is a point
$y\in S_{2d}$ such that $\SS_y$ is Cohen-Macaulay and of degree and dimension
the same as the secant variety of an elliptic normal curve of degree $2d$,
then for the general $y\in S_{2d}$, $\SS_y$ must coincide with
the secant variety $Sec(E)$, where $y\in E$. Thus in order to prove the
theorem,
we need only to exhibit one such point.

If $d>3$, consider $(y_0:\ldots:y_{2d-1}):=(1:1:0:\ldots:0)$. Then
for this choice of the parameter
$$M_d=\pmatrix{x_0&0&0&\cdots&0&x_{2d-1}\cr
x_1&x_2&0&\cdots&0&0\cr
0&x_3&x_4&\cdots&0&0\cr
&\ldots&&&\ldots&\cr
0&0&0&\cdots&x_{2d-3}&x_{2d-2}\cr}.$$
This matrix has the property that any of its $3\times 3$ minors is either zero
or consists of one monomial. Furthermore, if $d(i,j)=1$, then $x_i$ and $x_j$
appear in the same row or column, so no $3\times 3$ minor can be of the form
$x_ix_jx_k$ with $d(i,j)\le 1$. Conversely, if $x_ix_jx_k$ is a monomial
such that $d(i,j)\ge 2$, $d(i,k)\ge 2$, $d(j,k)\ge 2$, then
$x_i$, $x_j$ and $x_k$ appear in distinct rows and columns in the above matrix,
and thus there is a $3\times 3$ minor of $M_d$ of the form $x_ix_jx_k$.
Therefore the ideal $I_d$ generated by the  $3\times 3$ minors of
$M_d$ is precisely $I_{X(\partial C(2d,4))}$,
by Proposition 4.1. In particular, $I_d$ cuts out a
projectively Gorenstein scheme of dimension 3 and
degree ${1\over 2}(2d)(2d-3)$.
This is exactly what we had to show.
$\bullet$

In the case of the secant variety of a degree $5$ elliptic normal curve,
$Sec(E)$ is a well known quintic hypersurface in $\Pfour$, and its geometry
has been extensively studied in the literature beginning with [Seg].
We refer to [ADHPR1] and [ADHPR2] for details and for
further references to the literature.
For an explicit form of its equation see [Hu], p. 109, or
[ADHPR1], [ADHPR2]. We give in the following a slightly different
version of this quintic  equation (compare [ADHPR2],
Proposition 4.12, $i$).

\proclaim Theorem 5.3. Let $y=(y_0:y_1:y_2:y_3:y_4)\in E\subset\Pfour$
be a general point. Then
$$Sec(E)=\{\det (x_{3(i+j)}y_{3(i-j)})_{i,j\in\boldz_5}=0\}.$$

Proof. By Corollary 2.9, the $3\times 3$ minors of the matrix
$M'_2=(x_{3(i+j)}y_{3(i-j)})_{i,j\in\boldz/5\boldz}$ vanish along $E$,
and thus the quintic $Q=\det(x_{3(i+j)}y_{3(i-j)})$ is singular
with multiplicity at least 3 along
$E$. In particular, $\{Q=0\}$ contains every secant line to $E$.
Since on the other hand, $Sec(E)$ is a hypersurface of degree 5,
we deduce that $Sec(E)=\{Q=0\}$.
$\bullet$

\proclaim Theorem 5.4. Let $E\subseteq \P^{2d}$ be a Heisenberg invariant
elliptic normal curve of degree $2d+1$, $d\ge 3$. Then the
intersection $E\cap (\P^{d-1})^-$ is non-empty, and for
a point $y=(y_0:\ldots:y_{2d})\in E\cap (\P^{d-1})^-$ the matrix
$$M'_{d}=(x_{(d+1)(i+j)}y_{(d+1)(i-j)})_{i,j\in\boldz/(2d+1)\boldz}$$
is skew-symmetric, and the ideal $I_{2d+1}$ generated by the
$6\times 6$-pfaffians of $M'_{d}$ is the ideal of $Sec(E)$,
for the general curve $E$.

Proof. By [LB], Proposition 4.7.5, $E\cap (\P^{d-1})^-$ is non-empty. Now
a point $y=(y_0:\cdots:y_{2d})$ in $(\P^{d-1})^-$ satisfies
$y_i=-y_{-i}$, for all $i\in\boldz_d$, whence
$$\eqalign{(M'_{d})_{ij}&=x_{(d+1)(i+j)}y_{(d+1)(i-j)}\cr
&=-x_{(d+1)(i+j)}y_{(d+1)(j-i)}\cr
&=-(M'_{d})_{ji}.\cr}$$
Therefore $M'_{d}$ is skew-symmetric, and by Corollary 2.9,
the $4\times 4$-pfaffians of $M'_{d}$ vanish along the curve $E$.
In particular, the $6\times 6$-pfaffians of $M'_{d}$ are cubic
hypersurfaces which are singular along $E$, and hence contain
$Sec(E)$. By taking $y=(y_i)$, whose coordinates are
$$y_i=\cases{1&if $i=d$;\cr
-1&if $i=d+1$;\cr
0&otherwise,\cr}$$
and which is a point lying on the standard $(2d+1)$-gon $X(\Gamma_{2d+1})$,
we obtain
$$M'_{d}=\pmatrix{
0&x_{d+1}&0&0&\cdots&0&0&-x_d\cr
-x_{d+1}&0&x_{d+2}&0&\cdots&0&0&0\cr
0&-x_{d+2}&0&x_{d+3}&\cdots&0&0&0\cr
&\vdots&&&&&\vdots&\cr
0&0&0&0&\cdots&-x_{d-2}&0&x_{d-1}\cr
x_d&0&0&0&\cdots&0&-x_{d-1}&0\cr}$$
It is now easily seen that the $6\times 6$-pfaffians of the above
matrix consist of all monomials $x_ix_jx_k$ with $d(i,j)\ge 2$,
$d(i,k)\ge 2$, and $d(j,k)\ge 2$. Indeed, to obtain the monomial
$x_ix_jx_k$ with this restriction, we have to take the minor whose
rows (and columns) are indexed by the set
$$\{i+d,i+d+1,j+d,j+d+1,k+d,k+d+1\} \subseteq \{0,\ldots,2d\},$$
where all indices are computed mod $2d+1$. This yields
the minor (except for the case that $d\in \{i,j,k\}$, when
a slightly different configuration results, which may be handled similarly):
$$\pmatrix{
0&x_i&0&0&0&0\cr
-x_i&0&\alpha&0&0&0\cr
0&-\alpha&0&x_j&0&0\cr
0&0&-x_j&0&\beta&0\cr
0&0&0&-\beta&0&x_k\cr
0&0&0&0&-x_k&0\cr}$$
where the values $\alpha$ and $\beta$ are irrelevant.
The pfaffian of this matrix is $x_ix_jx_k$. Furthermore,  a minor
which does not contain three adjacent pairs of rows, must contain
a zero row and hence its pfaffian is identically zero.

Therefore, by Proposition 4.1, we deduce
that $I_{2d+1}=I_{X(\partial C(2d+1,4))}$. Moreover, using the same
argument as in the even case,
this again implies that for the general elliptic normal curve $E$ and
$y\in E\cap(\P^{d-1})^-$, the homogeneous ideal $I_{2d+1}$ cuts out
$Sec(E)$.
$\bullet$

Finally the results in \S 4 allow us also to say something about
minimal free resolutions. Part $a)$ below is  a rather
well-known fact (see e.g. [E], Exercise A2.22):

\proclaim Theorem 5.5.  Let $E\subseteq \P^{n}$ be a Heisenberg invariant
elliptic normal curve. Then
\item{$a)$} For $n\ge 4$ the homogeneous ideal $I_E$ has
a minimal free resolution of type
$$0\leftarrow I_{E}\leftarrow
R(-2)^{b_1}\leftarrow R(-3)^{b_2}\leftarrow\cdots
\leftarrow
R(-n+1)^{b_{n-2}}\leftarrow R(-n-1)\leftarrow 0,$$
where
$$b_i=i{{n}\choose {i+1}}-
{{n-1}\choose {i-1}},$$
for all $1\le i\le n-2$, and  where $R$ is the homogeneous
coordinate ring of $\P^n$.
\item{$b)$} $Sec(E)$ is projectively Gorenstein and its  homogeneous ideal
$I_{Sec(E)}$ has a minimal free resolution of type
$$0\leftarrow I_{Sec(E)}\leftarrow
R(-3)^{b_1}\leftarrow R(-4)^{b_2}\leftarrow\cdots
\leftarrow
R(-n+2)^{b_{n-4}}\leftarrow R(-n-1)\leftarrow 0,$$
if $E$ is general and $n\ge 6$, where
$$b_i={{n-1}\choose {i+2}}{{i+1}\choose 2}+
{{n-1}\choose i}{{n-2-i}\choose 2} -
{{n-3}\choose i}{{n-2}\choose 2},$$
for all $1\le i\le n-4$.\par

Proof. To prove $b)$ it is enough to check the claims in
a special case, for instance the secant variety of $X(\Gamma_{n+1})$
which in turn (as  seen in \S 4) may be identified with the face
variety $X(\partial C(n+1,4))$. As mentioned in \S4, its homogeneous
coordinate ring is  an
extremal compressed Gorenstein algebra, and thus our claims follow from
[Sch], Theorem B (see also [FL]),
which essentially asserts that the Betti numbers of such extremal rings
are uniquely determined by their Hilbert function.
Part $a)$ may  also be regarded as a corollary of this result since
homogeneous rings of elliptic curves are obviously
extremal.$\bullet$

\par\bigskip

\noindent
{\hd \S 6. Equations defining abelian surfaces and the structure
of $\A_n^{lev}$.}

The goal in this section is to understand the structure of the ideal
of a general abelian surface embedded via a polarization of type $(1,n)$,
for $n\ge 10$. In this case, the homogeneous ideal will be generated
by quadrics which are easy to write down.\par\bigskip

{\it $n$ even:} We write $n=2d$. In this case, the projectivization of
the negative eigenspace of $\iota$ is $\P^-\cong \P^{d-2}$.

\proclaim Lemma 6.1. \item{$a)$} Let $A$ be an abelian surface with a line
bundle $\L$ of type $(1,2d)$, $d\ge 2$, $\L$ being of
characteristic zero with respect to a given decomposition,
and let $\phi:A\rightarrow \P^{2d-1}$
be the morphism induced by a basis of canonical theta functions. Then
$\phi(A)\cap\P^-$ consists of four distinct points.
\item{$b)$} Let $S_{E,\tau}$ be a Heisenberg invariant translation scroll
in $\P^{2d-1}$, with $E$ an elliptic curve in $\P^{2d-1}$ of degree $2d$. Then
$S_{E,\tau}\cap\P^-$ consists also of four distinct points.

Proof: $a)$ Let $A=V/\Lambda$, and let $\lambda_1,\lambda_2,\mu_1,\mu_2$
be a symplectic basis for $\Lambda$ compatible with a given decomposition,
with $E(\lambda_1,\mu_1)=1$, $E(\lambda_2,\mu_2)=2d$, following the notation
of \S 1.
Then for a 2-torsion point $v\in A$, represented by $\lambda/2$,
$\lambda\in\Lambda$, we define
$$q_{\L}(v):=\exp({\pi i E(\lambda',\mu')}),$$
where $\lambda=\lambda'+\mu'$ is the decomposition of $\lambda$ with
respect to the chosen decomposition of $\Lambda$. It is then
easy to see from explicit calculations that
$$\eqalign{\#\{v\in A_2\mid q_{\L}(v)=+1\}=12,&\cr
\#\{v\in A_2\mid q_{\L}(v)=-1\}=4.&\cr}$$
By [LB], Proposition 4.7.2, for $v\in A_2$, if $D$ is a symmetric divisor, then
$$(-1)^{mult_v(D)-mult_0(D)}=q_L(v).$$
So if $D$ is an even symmetric divisor, that is $mult_0(D)$ is even,
we see that
$$\eqalign{\#\{v\in A_2\mid\hbox{$mult_v(D)$ even}\}=12,&\cr
\#\{v\in A_2\mid\hbox{$mult_v(D)$ odd}\}=4.&\cr}$$
The claim then follows, since $0\in A$ is mapped to $\P^+$.

$b)$ Let $E\subseteq \P^{2d-1}$ be a Heisenberg invariant elliptic normal
curve. Let $\langle x,y\rangle$ be a secant line of $E$ joining the points
$x,y\in E$. If $\langle x,y\rangle$ intersects $\P^-$, then
$\langle x,y \rangle$ and $\iota\langle x,y \rangle$ span at most a
$\Ptwo$. But a $\Ptwo$ cannot intersect $E$ in four points, so either $x$
or $y$ is fixed by $\iota$ or the whole line $\langle x,y \rangle$ is fixed
by $\iota$. In the former case, this tells us that either $x$ or $y$ is
in $(\P^d)^+$, since $E\cap \P^-=\emptyset$ by [LB], Corollary 4.7.6.
Thus two points of
$\langle x,y \rangle$ are fixed by $\iota$, and so the line
$\langle x,y \rangle$ is fixed by $\iota$ anyway. Hence
$\langle x,y \rangle\cap \P^-\not=\emptyset$ if and only if $x=-y$
in the group law on  $E$.
Now a translation scroll $S_{E,\tau}$ is defined as
$$\bigcup_{x\in E} \langle x, x+\tau \rangle,$$
so $x+\tau=-x$ if and only if $2x=-\tau$, and there are precisely four such
possible values of $x$.
$\bullet$

Both $\sigma^d$ and $\tau^d$ act on $\P^-$, and this defines a $\boldz_2
\times\boldz_2$ action on $\P^-$. Thus, by Lemma 6.1, if $A$ is a Heisenberg
invariant abelian surface or translation scroll in $\P^{2d-1}$, then
$A\cap\P^-$ consists of a $\boldz_2\times\boldz_2$ orbit on $\P^-$.

Using Proposition 1.3.1, we can identify a point of $\A_{2d}^{lev}$
with a Heisenberg invariant
abelian surface $A$ contained in $\P^{2d-1}$, if the restriction of the
universal line bundle is very ample on $A$. Based on this, we give the
following definition:

\proclaim Definition. For $d\ge 2$, we may define a rational map
$$\Theta_{2d}:\A_{2d}^{lev}\rto \P^-/\boldz_2\times\boldz_2$$
by mapping $A$ to the $\boldz_2\times\boldz_2$ orbit $A\cap\P^-$.

\proclaim Theorem 6.2.
\item{$a)$} For a general $\HHH_{2d}$-invariant, $(1,2d)$-polarized
 abelian surface $A\subseteq
\P^{2d-1}$, and a point $y\in A\cap \P^-$, the $4\times 4$-pfaffians of
the matrices
$$\cases{
M_5(x,y), M_5(x,\sigma^5(y)),M_5(x,\tau^5(y))&if $d=5$;\cr
M_d(x,y), M_d(x,\sigma^d(y))& if $d\ge 7$, $d$ odd;\cr
M_6(x,y), M_6(\sigma(x),y), M_6(\tau(x),y)& if $d=6$;\cr
M_d(x,y), M_d(\sigma(x),y)& if $d\ge 8$, $d$ even\cr}$$
generate the homogeneous ideal of
$A$.
\item{$b)$} For $d\ge 5$, $\Theta_{2d}$ is birational onto its image.

Proof: Let $\pi:\P^-\rightarrow \P^-/\boldz_2\times\boldz_2$ be the quotient
map, and let $Z_{2d}:=\pi^{-1}(\im\Theta_{2d})$. Let also $S_{E,\tau}$ be a
translation scroll. By Theorem 3.1, there exists a flat family $\A\rightarrow
\Delta$, $\A\subseteq\P_{\Delta}^{2d-1}$, of Heisenberg invariant surfaces
with $\A_0=S_{E,\tau}$, and $\A_t$ a non-singular $(1,2d)$-polarized
abelian surface for the
general $t\in\Delta$. Then it is clear that $S_{E,\tau}\cap\P^-
\subseteq\overline{Z_{2d}}$. We deduce that for any elliptic normal
curve, $Sec(E)\cap
\P^-\subseteq\overline{Z_{2d}}$, and so the same is true if instead of a
non-singular elliptic curve, we take $E$ to be the ``standard $2d$-gon''
and $Sec(E)$ to be its ``secant variety''. Thus, in particular, the point
$$p_0:=(0:1:1:0:\cdots:0:-1:-1)$$
lies in $\overline{Z_{2d}}$.
Consider now the family
$$\A\subseteq\P^{2d-1}\times\overline{Z_{2d}}$$
over $\overline{Z_{2d}}$ defined by the $4\times 4$-pfaffians of the matrices
$M_d(x,y)$, etc (where the choice is made depending on $d$ as
in the statement of the theorem).
Here, the $x$'s are coordinates on the first $\P^{2d-1}$ while the $y$'s
are coordinates for $\overline{Z_{2d}}\subseteq \P^{2d-1}$.
By Corollary 2.7, $\A_y$ contains all
abelian surfaces represented by the points of $\Theta_{2d}^{-1}(\pi(y))$,
for $y\in Z_{2d}$.

Now some easy combinatorics show that the $4\times 4$-pfaffians
of the matrices in the statement of the Theorem
generate the ideal $I_{X_{2d}^{1}}$, using the description of
this ideal given in Theorem 4.6. We prove this in the general
case for $d\ge 7$, and leave it to the reader to check the cases
$d=5$ and $d=6$. We first consider certain pfaffians of $M_d(x,y_0)$.
If $i$ is odd, $0\le i<2d$, then the $4\times 4$-submatrix consisting
of the rows and columns indexed by $${i-1\over 2}, {i-1\over 2}+1,
{i-1\over 2}+2,{i-1\over 2}+3,$$ modulo $d$, is the matrix
$$\pmatrix{
0&-x_i&-x_{i+1}&0\cr
x_i&0&-x_{i+2}&-x_{i+3}\cr
x_{i+1}&x_{i+2}&0&-x_{i+4}\cr
0&x_{i+3}&x_{i+4}&0\cr}$$
whose pfaffian is easily seen to be
$$x_ix_{i+4}-x_{i+1}x_{i+3}.$$
If $i$ and $j$ are both odd, and $d(i,j)\ge 6$, then the submatrix
consisting of the rows and columns
$${i-1\over 2}, {i+1\over 2}, {j-1\over 2}, {j+1\over 2}$$ is now
$$\pmatrix{
0&-x_i&0&-\alpha\cr
x_i&0&-\beta&0\cr
0&\beta&0&-x_j\cr
\alpha&0&x_j&0\cr}$$
where $\alpha$ and $\beta$ depend on the precise choice of $i$ and
$j$, but at least one of them is zero. The pfaffian of this matrix is
$x_ix_j$. Finally, if $i$ is odd and $j$ even, with $d(i,j)\ge 5$,
then the submatrix with rows
$${i-1\over 2}, {i+1\over 2}, {j\over 2}-1, {j\over 2}+1$$ is
$$\pmatrix{
0&-x_i&-\alpha&-\beta\cr
x_i&0&-\gamma&-\delta\cr
\alpha&\gamma&0&-x_j\cr
\beta&\delta&x_j&0\cr}$$
with $\alpha,\ldots,\delta$ depending on the precise choice of
$i$ and $j$, but with $\alpha\delta=\beta\gamma=0$; in particular
the pfaffian of this matrix is $x_ix_j$. In this way, we have
obtained the following sets of equations
$$\{x_{i+1}x_{i+3}-x_ix_{i+4}\mid\hbox{$i$ odd}\}\cup
\{x_ix_j\mid \hbox{$d(i,j)\ge 5$, $d(i,j)$ odd}\}\cup$$
$$\{x_ix_j\mid \hbox{$d(i,j)\ge 6$, $d(i,j)$ even, $i$, $j$ odd}\}.$$
To obtain the remaining equations, if $d$ is even, we  apply
$\sigma$ to these equations. If $d$ is odd, we apply a similar procedure
to the matrix $M_d(x,\sigma^d(p_0))$.

Recall now from Theorem 4.6 that
$X_{2d}^{1}$ is a Cohen-Macaulay
surface of degree
$4d$. Since for a general $y\in Z_{2d}$, $\A_y$ is either of dimension
$>2$, or is a surface of degree at least $4d$, we see that for
general $y$, $\dim\A_y=2$ and $\deg\A_y=4d$. Since furthermore
$\A_{p_0}=X_{2d}^{1}$ has the same Hilbert polynomial as a $(1,2d)$-polarized
abelian surface, there must exist an open neighborhood
$U\subseteq\overline{Z_{2d}}$ of $p_0$ such that $\A_U\subseteq
\P^{2d-1}_U$ is flat over $U$, and each smooth fibre of $\A_U\rightarrow U$
is an abelian surface. In particular, since the ideal of $\A_{p_0}$ is
generated by the pfaffians in question, the same is true for
the ideal of $\A_y$ for
general $y\in U$. Thus $a)$ follows, and $b)$ is now clear since
$\Theta_{2d}:\Theta^{-1}_{2d}(\pi(U\cap Z_{2d}))
\rightarrow \pi(U\cap Z_{2d})$ is an isomorphism.
$\bullet$
\par\bigskip

{\it $n$ odd:} Let $n=2d+1$. The projectivization of the negative
eigenspace of the involution $\iota$ acting on $\P^{2d}$ is
now $\P^-\cong\P^{d-1}$.
The group $\HHH_{2d+1}$ acts on $H^0(\O_{\P^{2d}}(2))$, and one
sees readily that
$H^0(\O_{\P^{2d}}(2))$ splits up into $d+1$ mutually isomorphic irreducible
representations of $\HHH_{2d+1}$. Let $R_d$ be the $(d+1)\times
(2d+1)$ matrix given by
$$(R_d)_{ij}=x_{j+i}x_{j-i},\quad 0\le i \le d,\ 0\le j\le 2d.$$
The rows of this matrix each span an irreducible subrepresentation
of
$H^0(\O_{\P^{2d}}(2))$
and yield the decomposition of
$H^0(\O_{\P^{2d}}(2))$
into $d+1$ mutually isomorphic representations.

Let $D_i\subseteq\P^-$ be the locus in $\P^-$ where $R_d$ has rank
$\le 2i$. Note that on $\P^-$, $x_0=0$ and $x_i=-x_{-i}$ for
$i\not=0$, so using coordinates $x_1,\ldots,x_d$ on $\P^-$, $R_n$ restricted
to $\P^-$ has the property that for $j\not =0$, the $j$-th column and
the $(2d+1-j)$-th column coincide, since
$$\eqalign{(R_d)_{i,2d+1-j}&=x_{(2d+1-j)+i}x_{(2d+1-j)-i}\cr
&=(-x_{j-i})(-x_{j+i})\cr
&=x_{j-i}x_{j+i}\cr
&=(R_d)_{ij}\cr}$$
on $\P^-$.
Also, the leftmost $(d+1)\times(d+1)$ block of $R_d$ when restricted to
$\P^-$ is anti-symmetric, since on $\P^-$:
$$\eqalign{(R_d)_{ij}&=x_{j+i}x_{j-i}\cr
&=x_{j+i}(-x_{i-j})\cr
&=-(R_d)_{ji}.\cr}$$
We denote by $T_d$ the restriction of this $(d+1)\times(d+1)$
block to $\P^-$. $D_i$ is then the locus where $T_d$ is rank $\le 2i$.

\proclaim Lemma 6.3. For a general $\HHH_{2d+1}$-invariant abelian
surface $A\subseteq \P^{2d}$, $d\ge 3$, we have $A\cap\P^-\subseteq D_2$ and
$A\cap\P^- \not\subseteq D_1$.

Proof: Note that $R_d$, up to transpose and permutations of rows and
columns, is a submatrix of the matrix $M'_{d}(x,x)$ of Corollary
2.8, which is rank at most 4 on $A$. Thus $A\cap\P^-\subseteq D_2$.
We need then to show that for general $A$, $A\cap\P^-\not\subseteq
D_1$. Let $S_{E,\tau}$ be a translation scroll. By Theorem
3.1, we can find a flat family $\A\rightarrow\Delta$ in $\P^{2d}_{\Delta}$
with $\A_0=S_{E,\tau}$, and such that $\A_t$ a non-singular abelian surface for
general $t\in\Delta$. Thus it is enough to show that
$S_{E,\tau}\cap\P^-\not\subseteq D_1$ for general $S_{E,\tau}$, or
equivalently, it is enough to show that $Sec(E)\cap\P^-\not\subseteq D_1$
for the general elliptic normal curve $E\subset\P^{2d}$.
By Theorem 3.2, it is then enough to check that
$Sec(E)\cap\P^-\not\subseteq D_1$ for $E=X(\Gamma_{2d+1})$,
the standard $(2d+1)$-gon. This
is clear since  the point
$$p_0=(0:1:1:0:\cdots:0:-1:-1)$$
lies in $Sec(E)\cap\P^-=X(\partial C(2d+1,4))\cap\P^-$, while the matrix
$$T_d(z_0)=\pmatrix{
0&1&1&0&\cdots&0&0&0\cr
-1&0&0&0&\cdots&0&0&0\cr
-1&0&0&0&\cdots&0&0&0\cr
0&0&0&0&\cdots&0&0&0\cr
&&\vdots&&&\vdots&\cr
0&0&0&0&\cdots&0&0&0\cr
0&0&0&0&\cdots&0&0&-1\cr
0&0&0&0&\cdots&0&1&0\cr},$$
has rank $4$.
$\bullet$

We may define now a morphism
$$s:D_2\backslash D_1\rightarrow Gr(d-3,d+1)$$
by taking $x\in D_2\backslash D_1$ to $\ker(R_d(x)^t)$. We think of the
points of
$Gr(d-3,d+1)$ as parametrizing $(d-3)(2d+1)$-dimensional
sub-$\HHH_{2d+1}^e$-representations of $H^0(\O_{\P^{2d}}(2))$: if
$N$ is the $(d+1)\times (d-3)$-matrix whose columns span a
$(d-3)$-dimensional subspace of ${\bf C}^{d+1}$, then the entries of
$R_d^t\cdot N$ span the corresponding $(d-3)(2d+1)$-dimensional
subrepresentation of $H^0(\O_{\P^{2d}}(2))$. In other words, the map $s$
takes a point $x\in D_2\backslash D_1$ to the largest subrepresentation of
quadrics in $H^0(\O_{\P^{2d}}(2))$ vanishing at $x$.

\proclaim Lemma 6.4. For a general $\HHH_{2d+1}$-invariant abelian
surface $A\subseteq\P^{2d}$, $d\ge 3$, the set
$s(A\cap\P^-\cap (D_2\backslash D_1))$ consists
of exactly one point $p\in Gr(d-3,d+1)$. If $V_p$ is the corresponding
$(d-3)$-dimensional subspace of ${\bf C}^{d+1}$, and $N_p$ is a matrix
whose columns span $V_p$, then the space of quadrics spanned by the
entries of $R_d^t\cdot N_p$ is $H^0(\I_A(2))$, and each column of
$R_d^t\cdot N_p$ spans a sub-$\HHH_{2d+1}$-representation of
$H^0(\I_A(2))$.

Proof: By Riemann-Roch,
$$\eqalign{h^0(\I_A(2))&\ge h^0(\O_{\P^{2d}})-h^0(\O_A(2))\cr
&= {2d+2 \choose 2} -4(2d+1)\cr
&=(d-3)(2d+1).\cr}$$
If $x\in A\cap\P^-$, $x\in D_2\backslash D_1$, and $p:=s(x)$, then the
entries of $R_d^t\cdot N_p$ span the largest sub-$\HHH_{2d+1}$-representation
of quadrics in $H^0(\O_{\P^{2d}}(2))$ vanishing at $x$. Call this
subrepresentation $I_2$. Now $\dim I_2=(d-3)(2d+1)$. Since
$H^0(\I_A(2))$ is also a sub-$\HHH_{2d+1}$-representation of
$H^0(\O_{\P^{2d}}(2))$
consisting of quadrics vanishing at $x$, we get
$H^0(\I_A(2))\subseteq I_2$. From the
dimension estimate above (or since the representations have weight 2,
cf. [La]), we must then have $H^0(\I_A(2))=I_2$. This is true
for each $x\in A\cap\P^-$, $x\in D_2\backslash D_1$, so we see that $s(x)$
is the point corresponding to the subrepresentation $H^0(\I_A(2))$
of $H^0(\O_{\P^{2d}}(2))$. In particular
$s(A\cap\P^-\cap (D_2\backslash D_1))$ consists of one point.
$\bullet$

\proclaim Definition. We define a rational map
$$\Theta_{2d+1}:\A_{2d+1}^{lev}\rto Gr(d-3,d+1)$$
by taking an abelian surface $A$ to
$s(A\cap\P^-\cap (D_2\backslash D_1))$, or equivalently, to the point of
$Gr(d-3,d+1)$ corresponding to the subrepresentation $H^0(\I_A(2))$ of
$H^0(\O_A(2))$.

\proclaim Theorem 6.5.
\item{$a)$} The homogeneous ideal of a general
$\HHH_{2d+1}$-invariant abelian surface
$A\subseteq\P^{2d}$ of type $(1,2d+1)$, $d\ge 5$, is
generated by quadrics.
\item{$b)$} $\Theta_{2d+1}$ is birational onto its image.

Proof: Let $Z:=s^{-1}(\im\Theta_{2d+1})\subseteq D_2\backslash D_1$, and let
$\overline{Z}$ be the closure of $Z$ in $D_2\backslash D_1$. Let $\A\subseteq
\P^{2d}_{\overline{Z}}$ be the family defined by the condition
that the ideal of $\A_z$, $z\in
\overline{Z}$, is the subrepresentation of $H^0(\O_{\P^{2d}}(2))$ determined by
$s(z)$, that is the subrepresentation of quadrics in
$H^0(\O_{\P^{2d}}(2))$ vanishing at $z$.
By Lemma 6.4, if $z\in s^{-1}(\Theta_{2d+1}(A))$, then $\A_z$ contains
the surface $A$.
Now exactly the same argument as in the proof of Lemma 6.3 shows
that the point
$$z_0=(0:1:1:0:\ldots:0:-1:-1)$$
is in $\overline{Z}$. Moreover, $\ker(R_d(z_0)^t)=\ker(T_d(z_0)^t)$, which
in turn is spanned
by the vectors $(0,1,-1,0,\ldots)$ and $e_4,\ldots,e_{d-1}$, where
$e_1,\ldots,e_{d+1}$ denote the standard basis of ${\bf C}^{d+1}$. Hence
$\A_{z_0}$ is defined by the ideal generated by the set of quadrics
$$\{x_{i+1}x_{i+3}-x_ix_{i+4}\mid i\in \boldz/(2d+1)\boldz\}
\cup \{x_ix_j\mid d(i,j)\ge 5\}.$$
By Theorem 4.6, the scheme $\A_{z_0}$ has the same Hilbert polynomial
as a non-singular $(1,2d+1)$-polarized abelian surface.
So an argument as in the proof
of Theorem 6.2 shows that there exists an open set $U\subseteq\overline{Z}$
such that $\A_U\rightarrow U$ is flat, and such that
each smooth fibre is an abelian surface;
the claims in $a)$ and $b)$ then follow.
$\bullet$
\par\medskip

\noindent
{\it Example 6.6.}
The results in Theorems 6.2 and 6.5 $a)$ do not hold for all
$(1,d)$-polarized abelian surfaces. More precisely, as the following
example shows, cubics may also be needed to generate the homogeneous ideal
(compare also [ADHPR1]):

Let $E:=\{x_0^3+x_1^3+x_2^3+\lambda x_0x_1x_2=0\}\subseteq\Ptwo=
\P(V)$, $\lambda\not\in\{\infty,-3,-3\epsilon_3,-3\epsilon_3^2\}$,
where $\epsilon_3$ is a third root of unity, be a smooth
cubic in the Hesse pencil (i.e., an elliptic normal curve in $\Ptwo$
embedded with canonical level structure). We also choose
as origin for the group law on $E$ an inflection
point, say $o_E=(0:1:-1)$.  Fix $d\in\boldz$, with $d\ge 3$ and $d\cong
2\mod 3$, and  let now $\{y_0, y_1,\ldots,y_{3d-1}\}$ be a
basis of canonical theta functions of $H^0(\O_E(3do_E))$. Recall that
the Heisenberg group $\HHH_{3d}$ acts on this basis via the Schr\"odinger
representation, namely
$$\sigma(y_i)=y_{i-1},\quad \tau(y_i)=\xi^{-i}y_i,\quad i\in\boldz_{3d},\
\xi:=\exp(2\pi i/3d).
$$
Since $[\sigma^d,\tau^d]=\xi^{-d^2}\cdot{\rm id}$, we may identify
$\HHH_3\subseteq SL(V)$ in its Schr\"odinger representation with
the subgroup of $\HHH_{3d}$ generated by $\sigma^d$ and $\tau^d$,
where we take as third root of unity $\epsilon_3:=\xi^{d^2}$.

Let now $E^\prime:=\{x_0^3+x_1^3+x_2^3+\lambda^\prime x_0x_1x_2=0\}
\subseteq\P(V^\ast)$, $\lambda^\prime\not\in
\{\infty,-3,-3\epsilon_3,-3\epsilon_3^2\}$,
with origin $o_{E^\prime}=(0:1:-1)\in\P(V^\ast)$,
and consider the diagonal action of $\HHH_3$ on the product
$E\times E^\prime\subseteq\P(V)\times\P(V^\ast)$, where the action
on the first factor is the one induced from $\HHH_{3d}$.
It is easy to see that the center of $\HHH_3$ acts
trivially on the line bundle $\M:=\O_E(3do_E)\boxtimes
\O_{E^\prime}(3o_{E^\prime})$, hence $\M$ descends to a line
bundle $\L:=\M/\boldz_3\times\boldz_3$ over the quotient
abelian surface $A$
defined via the 9-fold unramified cover
$$\pi:E\times E^\prime\to E\times E^\prime/\boldz_3\times\boldz_3=:A.$$
Since $\M$ is of type $(3,3d)$, and $\gcd(d,3)=1$, we deduce that $\L$
defines on $A$ a polarization of type $(1,d)$. By using Reider's criterion
it is also easily seen that $\L$ defines an embedding to $\P^{d-1}$,
whenever $d\ge 5$. A basis of sections for  $W=H^0(\L)$ is defined by the
following invariant sections of $\M$:
$$s_i=y_{3i}\boxtimes x_0+y_{3i+d}\boxtimes x_1+y_{3i+2d}\boxtimes x_2,\qquad
i\in\boldz_d,$$
where the indices of the $y's$ are taken modulo $3d$.
The subgroup of $\HHH_{3d}$
generated by $\sigma^3$ and $\tau^3$ acts naturally on $\M=\O_E(3do_E)\boxtimes
\O_{E^\prime}(3o_{E^\prime})$, where the action on the second factor is the
trivial one, and since it commutes with the above diagonal action
of $\HHH_3$ we deduce that this action descends to an action on the line
bundle $\L$, whence in particular on its sections. Moreover, since
$[\sigma^3,\tau^3]=\xi^{-9}\cdot{\rm id}$, we may identify this subgroup
with the Heisenberg group $\HHH_d\subseteq SL^\pm(W)$. It is also readily
checked
that $\HHH_d$ acts via the Schr\"odinger representation on the chosen basis
$\{s_i\mid i\in\boldz_d\}$ of $H^0(\L)$. In other words, we've picked this
way a canonical level structure on $(A,\L)$.

Let now $F$ and $F^\prime$ denote the images on $A$ of the elliptic curves
$E$ and $E^\prime$, respectively. Since $\deg \L_{|F^\prime}=3$,
all curves in the pencil $|F^\prime|$ are embedded via
$\L$ in $\P^{d-1}$  as plane cubic curves, and hence any quadric
hypersurface containing the abelian surface $\psi_{\L}(A)$
must also contain the threefold
traced out
by the planes spanned by the plane cubics in the pencil $|F^\prime|$.
In particular, we deduce that the homogeneous ideal of $\psi_{\L}(A)$
cannot be generated only by quadrics. It is also easy to check in this case
that for general choices of $E$ and $E^\prime$, the homogeneous ideal
of $\psi_{\L}(A)$ is generated by quadrics and cubics when $d\ge 7$.

The above results and rather extensive checks in examples have lead us
to formulate the following

\proclaim Conjectures.
\item{$a)$} The homogeneous ideal of an  embedded $(1,d)$-polarized abelian
surface
is generated by quadrics and cubics, for $d\ge 9$.
\item{$b)$} The general  $(1,d)$-polarized abelian surface, for $d\ge 10$, has
property $N_{\lbrack {d\over 2}\rbrack-4}$.\par

\par\bigskip\noindent
{\hd References}
\item{[ADHPR1]} Aure, A.B., Decker, W., Hulek, K., Popescu, S., Ranestad, K.,
``The Geometry of Bielliptic Surfaces in $\Pfour$'',
{\it Int. J. of Math.}, {\bf 4}, (1993) 873--902.
\item{[ADHPR2]} Aure, A.B., Decker, W., Hulek, K., Popescu, S., Ranestad, K.,
``Syzygies of Abelian and Bielliptic Surfaces in $\Pfour$'',
preprint alg-geom/9606013.
\item{[Au]} Aure, A.,
 ``Surfaces on quintic threefolds associated to the Horrocks-Mumford bundle'',
in {\it Lecture Notes in Math.}, {\bf 1399}, (1989), 1--9, Springer.
\item{[Ba]} Barth, W., ``Abelian surfaces with $(1,2)$ polarization'',
{\it Adv. Studies in Pure Math.}, {\bf 10}, Alg. Geom. Sendai, (1985) 41--84.
\item{[BaH]} Barth, W., Hulek, K., ``Projective models of Shioda
modular surfaces'', {\it Manuscripta Math.}, {\bf 50}, (1985) 73--132.
\item{[BE]} Bayer, D., Eisenbud, D., ``Graph Curves. With an appendix by
Sung Won Park'', {\it Adv. in Math.},{\bf 86}, No. 1, (1991) 1--40.
\item{[BH]} Bruns, W., Herzog, J., {\it Cohen-Macaulay Rings},  Cambridge
Studies in advanced mathematics, {\bf 39}, Cambridge University Press 1993.
\item{[BLvS]} Birkenhake, Ch., Lange, H., van Straten, D., ``Abelian
surfaces of type $(1,4)$'', {\it Math. Ann.}, {\bf 285}, (1989) 625--646.
\item{[BS]} Bayer, D., Stillman, M.,
``Macaulay: A system for computation in
        algebraic geometry and commutative algebra
Source and object code available for Unix and Macintosh
        computers''. Contact the authors, or download from
        {\bf math.harvard.edu} via anonymous ftp.
\item{[Cs]} Cs\'asz\'ar, A., ``A polyhedron without diagonals'',
{\it Acta Sci. Math. (Szeged)}, {\bf 13}, (1949) 140--142.
\item{[DHS]} Debarre, O., Hulek, K., and Spandaw, J.,
``Very ample linear systems on abelian varieties,'' {\it Math. Ann.}
{\bf 300}, (1994) 181-202.
\item{[DR]} Deligne, P., and Rapaport, M., ``Les sch\'emas de modules
de courbes elliptiques,'' in {\it Modular Functions of One Variable, II}
{\it Lecture Notes in Math.}, {\bf 349}, (1973), 143--316, Springer.
\item{[E]} Eisenbud, D., {\it Commutative Algebra with a View
Toward Algebraic Geometry}, Springer 1995.
\item{[EKS]} Eisenbud, D., Koh, J., Stillman, M., ``Determinantal
equations for curves of high degree'', {\it Amer. J. of Math.}, {\bf
110}, (1988) 513--539.
\item{[ES]} Eisenbud, D., Sturmfels, B., ``Binomial ideals'', {\it Duke
Math. J.}  {\bf 84}, (1996), no. 1, 1--45.
\item{[FL]} Fr\"oberg, R., Laksov, D., ``Compressed algebras'',
in ``Complete intersections, Acireale 1983'', {\it Lecture Notes in Math.}
{\bf 1092}, (1984), 121--151, Springer.
\item{[Gr]} Gr\"unbaum, B., {\it Convex Polytopes}, J. Wiley and Sons 1967.
\item{[GP]} Gross, M., Popescu, S., ``Calabi-Yau 3-folds and moduli of
abelian surfaces'', in preparation.
\item{[Ha]} Hartshorne, R., {\it  Algebraic Geometry,} Springer 1977.
\item{[Ho]} Hochster, M.,``Cohen-Macaulay rings, combinatorics and simplicial
complexes'', in ``Ring theory II'',  McDonald B.R., Morris, R. A. (eds),
{\it Lecture Notes in Pure and Appl. Math.}, {\bf 26}, M. Dekker 1977.
\item{[HM]} Horrocks, G., Mumford, D., `` A Rank 2 vector bundle on  $\Pfour$
with 15,000 symmetries'', {\it Topology}, {\bf 12}, (1973) 63--81.
\item{[HR]} Hochster, M., Roberts, J., ``The purity of the Frobenius
and local cohomology'', {\it Advances in Math.}, {\bf 21}, (1976),
no. 2, 117--172.
\item{[Hu]} Hulek, K., ``Projective geometry of elliptic curves'',
{\it Ast\'erisque}, {\bf 137} (1986).
\item{[HKW]} Hulek, K., Kahn, C., Weintraub, S., {\it Moduli Spaces of Abelian
Surfaces: Compactification, Degenerations, and Theta Functions},
Walter de Gruyter 1993.
\item{[Ke]} Kempf., G., ``Projective coordinate rings of abelian
varieties'', in {\it Algebraic Analysis, Geometry and Number Theory},
edited by J. Igusa, The John Hopkins Press 1989, 225--236.
\item{[Kl]} Klein, F., ``\"Uber transformationen siebenter Ordnung
der elliptischen Funktionen'' (1878/79), Abhandlung {\bf LXXXIV}, in
{\it Gesammelte Werke}, Bd. {\bf III}, Springer, Berlin 1924.
\item{[KlF]} Klein, F., Fricke, R., {\it Theorie der elliptischen
Modulfunktionen}, Bd. {\bf I}, Teubner, Leipzig 1890.
\item{[Ko]} Koizumi, S., ``Theta relations and projective normality of abelian
varieties'', {\it Am. J. Math.}, {\bf 98}, (1976) 865--889.
\item{[La]} Lazarsfeld, R., ``Projectivit\'e normale
des surfaces abeliennes'', (written by O. Debarre),
preprint Europroj {\bf 14}, Nice.
\item{[LB]} Lange, H., Birkenhake, Ch., {\it Complex abelian varieties},
Springer-Verlag 1992.
\item{[LN]} Lange, H., Narasimhan, M.S., ``Squares of ample line bundles on
abelian varieties'', {\it Expo. Math.}, {\bf 7}, (1989) 275--287.
\item{[Ma]} Manolache, N.,
``Syzygies of abelian surfaces embedded in ${\P^4({\bf C})}$'',
{\it J. Reine Angew. Math.}, {\bf 384}, (1988) 180--191.
\item{[Mo]} Moore, R.: ``Heisenberg invariant quintic 3-folds and
sections of the Horrocks-Mumford bundle'', preprint Canberra 1985.
\item{[Mu1]} Mumford, D., ``On the equations defining abelian varieties'',
{\it Inv. Math.}, {\bf 1}, (1966) 287--354.
\item{[Mu2]} Mumford, D., {\it Abelian varieties}, Oxford
University Press 1974.
\item{[Mu3]} Mumford, D., {\it Tata Lectures on Theta I},
 Progress in Math. {\bf 28}, Birkh\"auser, 1983.
\item{[MS]} Manolache, N., Schreyer, F.-O., private communication.
\item{[Oh]} Ohbuchi, A., ``A note on the normal generation of ample
line bundles on abelian varieties'', {\it Proc. Japan Acad.},
{\bf 64}, (1988) 119--120.
\item{[Ra]} Ranestad, K., private communication.
\item{[Rav]} Ravi, M.S., ``Determinantal equations for secant varieties
of curves'', {\it Comm. in Algebra}, {\bf 22}, (1994) 3103--3106.
\item{[Re]} Reisner, G.A., ``Cohen-Macaulay quotients of polynomial rings'',
{\it Adv. in Math.}, {\bf 21}, (1976) 30--49.
\item{[Ro]} Room, T.G., {\it The Geometry of the Determinantal Loci},
Cambridge Univ. Press 1938.
\item{[Sch]} Schenzel, P., ``\"Uber die freien Aufl\"osungen extremaler
Cohen-Macaulay-Ringe'', {\it J. of Algebra}, {\bf 64}, (1980) 93-101.
\item{[Seg]} Segre, C., ``Sull' Incidenza di rette e piani nello spazio
a quattro dimensioni'', {\it Rend. Palermo} {\bf II}, (1888) 42--52.
\item{[SR]} Semple, G., Roth, L., {\it Algebraic Geometry}, Chelsea 1937.
\item{[Sta]} Stanley, R., ``Combinatorics and Commutative Algebra, Second
edition'', Progress in Math. {\bf 41}, Birkh\"auser, 1996.
\item{[Ste]} Stevens, J., ``Degenerations of elliptic curves and cusps
singularities'', preprint (alg-geom/9512014).
\item{[Ve]} V\'elu, J., ``Courbes elliptiques munies d'un sous-groupe
$\boldz/n\boldz\times\mu_n$'', {\it Bull. Soc. Math. France},
M\'emoire {\bf 57}, (1978) 5--152.
\item{[Wh]} White, A., ``Graphs, Groups and Surfaces'',
{\it North-Holland Mathematics Studies}, {\bf 8}, North Holland 1984.
\end